\documentclass[reqno,11pt]{amsart}
\usepackage{amsmath, latexsym, amsfonts, amssymb, amsthm, amscd}
\usepackage{graphics,epsf,psfrag}
\numberwithin{equation}{section}

\newtheorem{theorem}{Theorem}[section]
\newtheorem{lemma}[theorem]{Lemma}

\newtheorem{cor}[theorem]{Corollary}
\newtheorem{rem}[theorem]{Remark}
\newtheorem{definition}[theorem]{Definition}

\renewcommand{\tilde}{\widetilde}

\newcommand{\cC}{{\ensuremath{\mathcal C}} }

\DeclareMathSymbol{\leqslant}{\mathalpha}{AMSa}{"36} 
\DeclareMathSymbol{\geqslant}{\mathalpha}{AMSa}{"3E} 
\DeclareMathSymbol{\eset}{\mathalpha}{AMSb}{"3F}     
\renewcommand{\leq}{\;\leqslant\;}                   
\renewcommand{\geq}{\;\geqslant\;}                   
\newcommand{\dd}{\,\text{\rm d}}             

\newcommand{\bbE}{{\ensuremath{\mathbb E}} }

\newcommand{\bbR}{{\ensuremath{\mathbb R}} }

\newcommand{\ga}{\alpha}
\newcommand{\gb}{\beta}
\newcommand{\gd}{\delta}
\newcommand{\gep}{\varepsilon}       

\newcommand{\gl}{\lambda}

\newcommand{\gs}{\sigma}

\makeatletter
\def\captionfont@{\footnotesize}
\def\captionheadfont@{\scshape}

\long\def\@makecaption#1#2{%
  \vspace{2mm}
  \setbox\@tempboxa\vbox{\color@setgroup
    \advance\hsize-6pc\noindent
    \captionfont@\captionheadfont@#1\@xp\@ifnotempty\@xp
        {\@cdr#2\@nil}{.\captionfont@\upshape\enspace#2}%
    \unskip\kern-6pc\par
    \global\setbox\@ne\lastbox\color@endgroup}%
  \ifhbox\@ne 
    \setbox\@ne\hbox{\unhbox\@ne\unskip\unskip\unpenalty\unkern}%
  \fi
  \ifdim\wd\@tempboxa=\z@ 
    \setbox\@ne\hbox to\columnwidth{\hss\kern-6pc\box\@ne\hss}%
  \else 
    \setbox\@ne\vbox{\unvbox\@tempboxa\parskip\z@skip
        \noindent\unhbox\@ne\advance\hsize-6pc\par}%
\fi
  \ifnum\@tempcnta<64 
    \addvspace\abovecaptionskip
    \moveright 3pc\box\@ne
  \else 
    \moveright 3pc\box\@ne
    \nobreak
    \vskip\belowcaptionskip
  \fi
\relax
}
\makeatother
\def\writefig#1 #2 #3 {\rlap{\kern #1 truecm
\raise #2 truecm \hbox{#3}}}

\newcommand{\dist}{\text{dist}}

\begin{document}

\title[Phase Reduction in the Escape Problem for Systems close to Reversibility]{Phase Reduction in the Noise Induced Escape Problem for Systems close to Reversibility}

\author{Christophe Poquet}
\address{
  Universit{\'e} Paris Diderot (Paris 7) and Laboratoire de Probabilit{\'e}s et Mod\`eles Al\'eatoires (CNRS),
U.F.R.                Math\'ematiques, Case 7012 
             75205 Paris Cedex 13, France
}

\begin{abstract}
We consider $n$-dimensional deterministic flows obtained by perturbing a gradient flow. We assume that the gradient flow admits a stable curve of stationary points, and thus if the perturbation is not too large the perturbed flow also admits an attracting curve. We show that the noise induced escape problem from a stable fixed point of this curve can be reduced to a one-dimensional problem: we can approximate the associated quasipotential by the one associated to the restricted dynamics on the stable curve. The error of this approximation is given in terms of the size of the perturbation.
\end{abstract}

\date{\today}

\maketitle

\begin{section}{Introduction}

\begin{subsection}{Phase reduction and escape problem}

For dynamical systems with an attracting limit cycle, the phase reduction method consists in simplifying the system by projecting the dynamics on the limit cycle,
and neglecting the distance between the trajectory and the limit cycle \cite{cf:Kur}. Such an approximation allows to reduce the dynamics to a one dimensional
self-contained equation satisfied by the phase parameterizing the limit cycle. Such a reduction is widely used in the context of noisy oscillators (see \cite{cf:GTNE,cf:TNE,cf:YA} and references therein).

The aim of this paper is to show that the phase reduction can be made in a rigorous way for the escape problem for a class of system close to reversibility.
For a smooth dynamical system 
\begin{equation}
\dd X_t \, =\, F[X_t] \dd t\, ,
\end{equation}
where $X_t\in\bbR^n$ (we use the notation $f[\cdot]$ for functions with domain $\bbR^n$), including a stable fixed point $A$ with basin of attraction $D$, the escape problem is the study of the metastable behavior of $A$ under a small noisy perturbation
\begin{equation}
\label{eq:LD generic}
\dd X_t\, =\,  F[X_t] \dd t +\sqrt{\gep}\dd B_t\, ,
\end{equation}
where $B_t$ is a Brownian motion in $\bbR^n$. The natural questions arising are where, when and how do trajectories of \eqref{eq:LD generic} escape from $D$.
This problem has been much studied in the literature. The fundamental reference is of course \cite{cf:FreidWentz}, where it is shown that these questions are related to the large deviation behavior
of \eqref{eq:LD generic}, and more precisely to the corresponding ``quasipotential''. For two points $P_1$ and $P_2$ of $\bbR^n$, the quasipotential $W(P_1,P_2)$ is defined by
\begin{equation}
 W(P_1,P_2)\, =\, \inf\left\{I_{T}^{P_1}(Y):\,Y\in C([T,0],\bbR^n),\, T<0,\,Y_{T}=P_1,\, Y_{0}=P_2\right\}\, , 
\end{equation}
where $I$ is the large deviation rate function of \eqref{eq:LD generic}, that is
\begin{equation}
 \label{eq:I general}
 I_{T}^x(Y)\, =\,\left\{ 
 \begin{array}{cl}
 \frac12\int_{T}^0 \left\Vert\dot Y_t-F[Y_t]  \right\Vert^2\dd t & \text{if } Y \text{ is absolutely continuous}\\ & \text{and } Y_{T}= x\, , \\
   +\infty & \text{otherwise}\, ,
 \end{array}
 \right.
\end{equation}
where $\Vert\cdot\Vert$ is the norm associated to the canonical scalar product $\langle \cdot, \cdot\rangle$ of $\bbR^n$.
For a compact neighborhood $K$ of $A$ with smooth boundaries included in $D$ (and thus attracted to $A$), it is shown in \cite{cf:FreidWentz} that the escape from $K$
will take place with probability tending to $1$ as $\gep\rightarrow 0$ very close to the points $B$ of the boundary of $K$ satisfying
$W(A,B)=\inf_{E\in \partial K}W(A,E)$. By a compactness argument and since $W$ is continuous \cite{cf:FreidWentz}, there exists at least one point $B$ satisfying this property. Moreover for each starting point $x\in K$, the exit time $\tau^\gep$ satisfies
\begin{equation}
\label{eq:limit exit time}
 \lim_{\gep\rightarrow 0} \gep \log \bbE_x \tau^\gep \, =\, W(A,B)\, . 
\end{equation}
Further work has been made to weaken the hypothesis on $K$, to study the escape at saddle points and the switching between the basins of
attraction of several stable fixed points (e.g. \cite{cf:Daybound,cf:GOV,cf:OV}). The escape time, after renormalization, is in fact
asymptotically exponential (e.g. \cite{cf:Dayexpo,cf:MOS,cf:OV}).

When \eqref{eq:LD generic} is reversible ($F=-\nabla V$ with $V$ smooth), the quasipotential is proportional to the potential driving the dynamics: if $B$ is in the basin of attraction of $A$, then
\begin{equation}
\label{eq:quas rev}
 W(A,B)\, =\, 2(V(B)-V(A))\, .
\end{equation}
In this case analytic approaches (in particular potential theory, see \cite{cf:BEGK,cf:Berglund2010}) show that the factor preceding $e^{W(A,B)/\gep}$ in \eqref{eq:limit exit time} satisfies the Eyring-Kramer's law \cite{cf:Eyring,cf:Kramer}.
We point out that in the one-dimensional case, since the escape problem only depends on the value of $F$ on a bounded domain, we can always consider that
the dynamics is driven by a gradient flow. 

\medskip

The purpose of this paper is to show that for a class of systems close to reversibility and containing an attracting curve $M$, the escape from a stable fixed point
$A$ located on $M$ occurs close to $M$, and that the quasipotential at the escape point can be approximated by the one corresponding to the dynamics constrained to $M$. We point out that the closeness to reversibility is a central point in our work. For what may happen far from reversibility, see for example \cite{cf:MS}: the most probable trajectories may go far away from the attracting curve, and the quasipotential cannot be reduced anymore.

In principle the results we prove here should also be true in infinite dimension, and this generalization would be particularly relevant (see \cite{cf:SellYou,cf:GPPP} for systems for which the infinite dimensional result would be of great interest).

\end{subsection}

\medskip

\begin{subsection}{Mathematical set-up and main result}

We will consider dynamical systems of the type
\begin{equation}
\label{eq:perturbed1}
 \dd X_t\, =\, (-\nabla V[X_t]+\gd G[X_t])\dd t +\sqrt{\gep}\dd B_t\, ,
\end{equation}
where $\gd$ is meant to be small, $G\in C^2(\bbR^n,\bbR^n)$ and $V\in C^4(\bbR^n,\bbR^n)$. 
The rate function associated to \eqref{eq:perturbed1} is
\begin{equation}
 \label{eq:I}
 I_{\gd,T}^x(Y)\, =\,\left\{ 
 \begin{array}{ll}
 \frac12\int_{T}^0 \Vert\dot Y_t+ \nabla V[Y_t]-\gd G[Y_t]  \Vert^2\dd t & \text{if } Y \text{ is absolutely continuous}\\ & \text{and } Y_{T}= x\, , \\
  +\infty & \text{otherwise}\, ,
 \end{array}
 \right.
\end{equation}
and we denote $W_\gd$ the associated quasipotential.

We suppose that the unperturbed deterministic dynamical system 
\begin{equation}
 \dd X_t\, =\, -\nabla V[X_t]\, ,
\end{equation}
contains a stable compact one-dimensional manifold of stationary solutions. More precisely, we suppose that there exists a curve $M$
($C^3$ by Local Inverse Theorem, since $V$ is $C^4$) without crossings such that for all $X\in M$ we have 
\begin{equation}
\label{eq:zero grad}
 \nabla V[X]\, =\, 0\, .
\end{equation}
For convenience we take $V\equiv 0$ on $M$. Moreover we will suppose the existence of a spectral gap for the linearized evolution in the neighborhood of $M$: we suppose that if $v$ is a tangent vector for $M$
at the point
$X$ and $w$ a vector orthogonal to $v$, then for $H[X]$ the Hessian matrix of $V$ at the
point $X$
\begin{equation}
\label{eq:zero hess}
 H[X]v\, =\, 0
\end{equation}
and there exists a positive constant $\gl$ (independant from the vector $w$) such that
\begin{equation}
\label{eq:spect gap hess}
 \langle H[X]w,w\rangle\, \geq\,  \gl\Vert w\Vert^2\, .
\end{equation}
With these hypothesis, $M$ is a normally hyperbolic manifold, which is a structure stable under small perturbations (see \cite{cf:Wig,cf:HPS,cf:SellYou}).
The perturbed deterministic dynamical system 
\begin{equation}
\label{eq:perturbed1 determ}
 \dd X_t\, =\, (-\nabla V[X_t]+\gd G[X_t])\dd t 
\end{equation}
thus also contains a stable normally hyperbolic curve $M^\gd$. However this new stable invariant manifold in general won't be a manifold of stationary solutions.
Moreover $M^\gd$ is located at distance of order $\gd$ from $M$ (more details will be given in Section \ref{subsec hyp}).

We consider a parametrization $\{q_\gd(\varphi),\, \varphi\in \bbR/L_\gd\bbR\}$ of $M^\gd$ satisfying $\Vert q'_\gd(\varphi)\Vert=1$ for all $\varphi$. $L_\gd$ denotes the length of the curve.
Since the flow \eqref{eq:perturbed1 determ} is tangent to $M^\gd$, for a trajectory $Y^\gd$ staying in $M^\gd$ (that is of the form $q_\gd(\varphi_t)$) linking two point $A^\gd$ and $B^\gd$ of $M^\gd$, the rate function is reduced to
\begin{equation}
\label{eq:I reduced}
 I_{\gd,T}^{A^\gd}(Y^\gd) \,=\, \int_{T}^0\Big|\dot\varphi^\gd_t-\Big\langle - \nabla V[q_\gd(\varphi_t)]+\gd G[q_\gd(\varphi_t)],q'_\gd(\varphi_t)\Big\rangle\Big|^2\dd t\, .
\end{equation}
This functional coincides with the large deviation rate function one obtains by considering the one-dimensional diffusion
\begin{equation}
\label{eq:1 dim reduction order gd}
 \dd \varphi^\gd_t\, =\, b_\gd\left(\varphi^\gd_t\right)\dd t+\sqrt{\gep}\dd B^1_t\, ,
\end{equation}
where $B^1$ is a one-dimensional Brownian motion and
\begin{equation}
\label{eq:def b gd}
 b_\gd(\varphi)\, :=\, \Big\langle -\nabla V[q_\gd(\varphi)]+\gd G[q_\gd(\varphi)],q'_\gd(\varphi)\Big\rangle\, .
\end{equation}
We denote $W^{red}_\gd(\varphi^1,\varphi^2)$ the associated quasipotential, i.e for all $\gd>0$, $\varphi_1\in \bbR/L_\gd\bbR$ and $\varphi_2\in \bbR/L_\gd\bbR$, $W^{red}_\gd$ is defined as follows:
\begin{multline}
 W_\gd^{red}(\varphi_1,\varphi_2)\, =\, \inf\bigg\{\int_{T}^0|\dot\varphi_t-b_\gd(\varphi_t)|^2\dd t:\,\varphi\in C([T,0],\bbR/L_\gd\bbR)\text{ and absolutely}\\ \text{continuous},\, T<0,\,\varphi_{T}=\varphi_1,\, \varphi_{0}=\varphi_2\bigg\}\, .
\end{multline}
Since $W^{red}_\gd$ is the infimum of the rate function taken on the subset made of the trajectories staying in $M^\gd$, we have the immediate bound
\begin{equation}
\label{ineq:comparaison W Wred}
 W_\gd(q_\gd(\varphi_1),q_\gd(\varphi_2))\, \leq \, W_\gd^{red}(\varphi_1,\varphi_2)\, .
\end{equation}
$b_\gd$ is in a certain sense a smooth perturbation of the function
\begin{equation}
\label{eq:def b}
b(\theta)\, =\, \gd\langle G[q(\theta)],q'(\theta)\rangle\, ,
\end{equation}
where $q$ is a parametrization of $M$ defined on $\bbR/L\bbR$ (where $L$ is the length of $M$) and satisfying $\Vert q'(\theta)\Vert=1$, so $b$ characterizes the perturbed dynamics projected on $M$.
The smoothness of the perturbation ensures that the dynamics induced by $b_\gd$ on $M^\gd$ will be conjugated to the one induced by $b$ on $M$ (and thus have the same properties). We consider the case when there exists a stable fixed
point $\theta^0$ for $b$ such that $b'(\theta^0)<0$, and such that the interval $[\theta^0-\Delta_1,\theta^0+\Delta_2]$ is included in the basin of attraction of $\theta^0$.
Then for $\gd$ small enough (see Lemma \ref{lem:tilde b} and the associated discussion), there exists a phase $\varphi_{A^\gd}$
(corresponding to a point $A^\gd=q_\gd(\varphi_{A^\gd})\in M^\gd$) also stable for $b_\gd$ and such that $[\varphi_{A^\gd}-\Delta_1,\varphi_{A^\gd}+\Delta_2]$
is included in its basin of attraction. Moreover, since $M^\gd$ is stable, $A^\gd$ is a stable fixed point also for \eqref{eq:perturbed1 determ}.

For each $Z$ close enough to $M^\gd$ there exists a unique $q_\gd(\varphi)$ such that $\Vert Z-q_\gd(\varphi)\Vert=\dist(Z,M^\gd)$. We denote by $p_\gd(Z):=\varphi$ the phase
of this projection. We will show that this phase gives the main contribution of the quasipotential associated to \eqref{eq:perturbed1} in the neighborhood of $M^\gd$.
More precisely, we define the tube
\begin{equation}
 \label{eq:def tube U}
  U^\gd\, =\, \{Z\in\bbR^n,\, \dist(Z,M^\gd)\leq C_0\gd^{1/2},\, p_\gd(Z)\in [\varphi_{A^\gd}-\Delta_1,\varphi_{A^\gd}+\Delta_2]\}\, ,
 \end{equation}
 depending on a constant $C_0$. We study the minimum of the quasipotential $W_\gd(A^\gd,\cdot)$ on the boundary $\partial U^\gd$ (recall that it is achieved), and the location of the points realizing this minimum. Since the ``length'' of such a tube may be of order $1$, whereas its ``slice'' is of order $\gd^{1/2}$, a trajectory exiting the tube at a point $B^\gd$ satisfying either
 $p_\gd(B^\gd)=\varphi_{A^\gd}-\Delta_1$ or $p_\gd(B^\gd)=\varphi_{A^\gd}+\Delta_2$ stays very close to $M^\gd$. We show that the most probable paths verify this property, and that the quasipotential $W_\gd$ can be well approximated by the reduced one, that is $W^{red}_\gd$.
 
\medskip

\begin{theorem}
\label{th:exp w}
There exist $\gd_0$ and a constant $C_0$ such that for all $\gd\leq\gd_0$, for each $B^\gd\in \partial U^\gd$ satisfying
\begin{equation}
\label{eq:min Wgd Bgd th}
 W_\gd\left(A^\gd,B^\gd\right)\, =\, \inf_{E\in \partial U^\gd}W_\gd\left(A^\gd,E\right)\, ,
\end{equation}
if we denote $\varphi_{B^\gd}:=p_\gd(B^\gd)$ we have either $\varphi_{B^\gd} =\varphi_{A^\gd}-\Delta_1$ or $\varphi_{B^\gd}=\varphi_{A^\gd}+\Delta_2$, and
\begin{equation}
\label{eq:th approx W}
W_\gd\left(A^\gd,B^\gd\right)\, =\,W^{red}_\gd\Big(\varphi_{A^\gd},\varphi_{B^\gd}\Big)+O\left(\gd^3|\log\gd|^3\right)\, .
\end{equation}
\end{theorem}

\medskip

This theorem proves that the quasipotential can be well approximated for the points satisfying the minimum of the quasipotential $W_\gd(A^\gd,\cdot)$ in the boundary of tube $U^\gd$. It is quite natural to think that this approximation is also possible for the points lying on the attracting curve $M^\gd$ of \eqref{eq:perturbed1 determ} (that is the points $B^\gd$ of the type $B^\gd=q_\gd(\varphi^\gd)$, but not necessarily satisfying \eqref{eq:min Wgd Bgd th}). This is the purpose of the following Corollary, obtained by carrying out a slight modification of the proof of Theorem \ref{th:exp w}.

\begin{cor}
\label{cor:w}
 There exist $\gd_0$ and a constant $C_0$ such that for all $\gd\leq\gd_0$ and for each $\varphi^\gd\in [\varphi_{A^\gd}-\Delta_1,\varphi_{A^\gd}+\Delta_2]$ we have
 \begin{equation}
W_\gd\left(A^\gd,q_\gd\left(\varphi^\gd\right)\right)\, =\,W^{red}_\gd\left(\varphi_{A^\gd},\varphi^\gd\right)+O\left(\gd^3|\log\gd|^3\right)\, .
\end{equation}
\end{cor}

\medskip

These results are obtained by quantitative estimates on the most probable paths. To understand why these paths stay at a distance of order $\gd^{1/2}$ from $M^\gd$ (or equivalently at distance $\gd^{1/2}$ from $M$, since $M^\gd$ is located at distance $\gd$ from $M$), remark that for a point $Z$ in the neighborhood of $M$, \eqref{eq:zero grad} and \eqref{eq:spect gap hess} imply that $V[Z]$ is equivalent up to a constant factor to $\dist(Z,M)^2$, where ``$\dist$'' denotes the distance associated to the norm $\Vert \cdot\Vert$. Since $A^\gd\in M^\gd$ and thus $\dist(A^\gd,M^\gd)=O(\gd^2)$, the contribution  to the quasipotential of the reversible part of the dynamics (see \eqref{eq:quas rev}) for such a point $Z$ is $V[Z]-V[A^\gd]=V[Z]+O(\gd^2)$. On the other hand, the fact that $b_\gd$ is of order $\gd$ shows that leaving $U^\gd$ following the curve $M^\gd$
has a cost of order $\gd$. This suggests that reaching a point located at a distance larger than $\gd^{1/2}$ is more expensive than following $M^\gd$. This idea is used in particular in the proof of Lemma \ref{lem:exit}.

\end{subsection}

\end{section}

\medskip

\begin{section}{Preliminary results of geometrical nature}\label{sec geom}

\begin{subsection}{Projection and local coordinates}\label{subsec proj}

We first give more details about the orthogonal projection on smooth curves. We are here in a particular case, since the manifold we want to project on is one-dimensional, and the
topology is induced by a scalar product. For the existence in more general cases, based on the Local Inverse Theorem, we refer for example to \cite{cf:Wig}. We will denote $\dist$ the distance associated to the norm $\Vert .\Vert$.

\begin{lemma}
\label{lem:p}
 Let $\cC$ be a $1$-dimensional $C^r$ manifold of $\bbR^n$ ($r\geq 2$). Let $s\mapsto g(s)$ be a $C^r$ parametrization of $\cC$ satisfying $\Vert g'(s)\Vert=1$.
 Then there exists a neighborhood $N$ of $\cC$ such that for all $Y$ in $N$ there exists a unique $s=p_g(Y)$ such that
 \begin{equation}
  \Vert Y-g(s)\Vert\, =\, \dist(Y,\cC)\, .
 \end{equation}
Moreover, for $s:=p_g(Y)$,
\begin{equation}
\label{eq:orthogonality proj}
 \langle Y-g(s),g'(s)\rangle\, =\, 0\, ,
\end{equation}
the mapping $Y\mapsto p_g(Y)$ is $C^{r-1}$, and for all $\gb\in \bbR^n$
\begin{equation}
 Dp_g[Y] \gb\, =\, \frac{1}{1-\langle Y-g(s), g''(s)\rangle}\langle g'(s),\gb\rangle \, .
\end{equation}
\end{lemma}

\medskip

\begin{proof}
The uniqueness of the projection for a sufficiently small neighborhood is ensured by the smoothness of $\cC$. \eqref{eq:orthogonality proj} is obtained by simply taking the derivative of
$\Vert Y-g(u)\Vert^2$ with respect to $u$ and the Implicit Function Theorem and \eqref{eq:orthogonality proj} imply that $p_g$ is $C^{r-1}$.
Let $h\in \bbR^n$ such that $\langle h, g'(\theta)\rangle=0$. Then it is clear that if $h$ is small enough such that the projection is well defined, $p_g(Y)=\theta$ for $Y=g(\theta)+h$. For a small perturbation $g(\theta)+h+\gb$, we are looking for the real $\ga$ satisfying
\begin{equation}
\label{eq:cond ga}
 \langle g(\theta)+h+\gb-g(\theta+\ga),g'(\theta+\ga)\rangle\, =\, 0\, .
\end{equation}
Since $p_g$ is $C^{r-1}$, we already know that $\ga=O(\Vert \gb\Vert)$. Now a first order expansion of \eqref{eq:cond ga} with respect to $\ga$ gives
\begin{equation}
\langle  -\ga g'(\theta) +h +\gb +O(\gb^2),g'(\theta)+\ga g''(\theta) +O(\gb^2)\rangle \, =\, 0\, ,
\end{equation}
which, since $\langle h, g'(\theta)\rangle =0$ and $\Vert g'(\theta)\Vert =1$ (which implies also $\langle g''(\theta),g'(\theta)\rangle =0$), leads to
\begin{equation}
 \ga(-1+\langle h, g''(\theta)\rangle)+\langle g'(\theta),\gb\rangle +O(\gb^2)\, =\, 0\, .
\end{equation}
\end{proof}

\medskip

In Theorem \ref{th:exp w} and in the rest of the paper, we consider a parametrization of $M$ (respectively of $M^\gd$) $\theta \mapsto q(\theta)$ for $\theta\in \bbR/L\bbR$
(respectively $\varphi\mapsto q_\gd(\varphi)$ for $\varphi\in\bbR/L_\gd\bbR$) satisfying $\Vert q'(\theta)\Vert=1$ (respectively $\Vert q_\gd'(\varphi)\Vert=1$) and we use the notations
\begin{align}
 p_\gd\, :=\, p_{q_\gd}\\
 p\, :=\, p_q\, .
\end{align}

\medskip

We stress out that the size of the neighborhood of a curve $\cC$ where the projection is defined depends continuously on its curvature and the sizes of its bottlenecks (which quantify in particular the non-crossing property of the curve). As we will see in Theorem \ref{lem:phi Mgd}, for the family of curves $M^\gd$ these quantities have continuous variations of order $\gd$. So if the projection $p$ is defined in a $\gep$-neighborhood of $M$ , this ensures the existence of the projections $p_\gd$ on a $(\gep+O(\gd))$-neighborhood of $M^\gd$ ($\gep$ fixed with respect to $\gd$), and in particular at distance $\gd^{1/2}$ from $M^\gd$ for $\gd$ small enough.

\end{subsection}

\medskip

\begin{subsection}{Stable Normally Hyperbolic Manifolds}

We now quickly review the notion of Stable Normally Hyperbolic manifolds (SNHM) (see \cite{cf:Wig} for more details). SNHMs are invariant manifolds, linearly stable, and such that the attraction they apply on their neighborhood is stronger than their inner dynamics. Consider a $C^r$ flow on $\bbR^n$
\begin{equation}
\label{eq:gen flow}
 \dot X\, =\, F(X)
\end{equation}
and suppose that it admits a compact invariant manifold $M$. Define for each $Q\in M$ its tangent space $T_Q$, its normal space $N_Q$ and the corresponding orthogonal projections $P^T_Q$ and $P^N_Q$.
To each initial condition $Q$ on $M$ we associate the linearized evolution semi-group $\Phi(Q,t)$ defined by
\begin{equation}
 \Phi(Q,0)u\, =\, u
\end{equation}
for all $u\in \bbR^n$ and
\begin{equation}
 \partial_t \Phi(Q,t)\, =\, DF(Q_t)\Phi(Q,t)
\end{equation}
where $Q_t$ is the trajectory of \eqref{eq:gen flow} with initial condition $Q$, and thus a trajectory staying on $M$. 

\begin{definition} For all $Q\in M$, we define the generalized Lyapunov-type numbers
\begin{equation}
 \nu(Q)\, :=\, \inf\left\{a:\, \left( \frac{\Vert w\Vert}{\Vert P^N_{Q_{t}}\Phi(Q,t)w\Vert} \right)\middle/ a^{-t}\rightarrow 0\quad\text{as}\quad t\downarrow -\infty\quad \forall w\in N_{Q}  \right\}
\end{equation}
and when $\nu(Q)<1$
\begin{multline}
 \gs(Q_0)\, :=\, \inf\Bigg\{b:\, \frac{\Vert w\Vert^b/\Vert v\Vert}{\Vert P^N_{Q_{t}}\Phi(Q,t)w\Vert^b/ \Vert P^T_{Q_{t}}\Phi(Q,t)v\Vert }\rightarrow 0 \\ \text{as}\quad t\downarrow -\infty\quad \forall v\in T_{Q},\, w\in N_{Q}  \Bigg\} \, .
\end{multline}

\end{definition}

The number $\nu$ characterizes the linear stability of $M$, and $\gs$ compares the normal and tangential linear evolution in the neighborhood of $M$. $\nu$ and $\gs$ are $C^r$ functions (see \cite{cf:Wig}), so they are bounded on the compact $M$, and attain their suprema $\bar\nu(M)$ and $\bar\gs(M)$ on $M$.

\begin{definition}
 $M$ is called a Stable Normally Hyperbolic Manifold if $\bar\nu(M)<1$ and $\bar\gs(M)<1$.
\end{definition}

It is clear that in our specific problem, the curve $M$ is a SNHM, since \eqref{eq:zero hess} and \eqref{eq:spect gap hess} imply $\bar\nu(M)\leq e^{-\gl}$ and $\bar\gs(M)=0$.
\end{subsection}

\medskip

\begin{subsection}{Persistence of hyperbolic manifolds}\label{subsec hyp}

We now formulate the persistence result of our 1-dimensional manifold $M$ under perturbation.  We refer to \cite{cf:fen,cf:Wig} for the general proof of persistence in the finite-dimension case. For more general cases (infinite dimensions), see for example \cite{cf:HPS,cf:SellYou}. Recall that $\theta\mapsto q(\theta)$ is a parametrization of $M$ satisfying $\Vert q'(\theta)\Vert =1$.

\medskip

\begin{theorem}
\label{lem:phi Mgd}
 If $G$ is $C^2$, then for all $\gd$ small enough, there exists a $C^2$ mapping $\theta \mapsto \phi_\gd(\theta)$ satisfying
\begin{equation}
\label{eq:orth phi}
 \langle \phi_\gd(\theta), q'(\theta)\rangle \, =\, 0\, ,
\end{equation}
\begin{equation}
 \sup_{\theta\in\bbR/L\bbR} \{\Vert\phi_\gd(\theta)\Vert,\, \Vert\phi'_\gd(\theta)\Vert,\, \Vert\phi''_\gd(\theta)\Vert\}\, =\, O(\gd)\, ,
\end{equation}
and such that
\begin{equation}
M^\gd\, =\, \{q(\theta) +\phi_\gd(\theta),\, \theta \in \bbR/L\bbR\}
\end{equation}
is a stable normally hyperbolic manifold for \eqref{eq:perturbed1 determ}.
\end{theorem}

\medskip

This result implies in particular that $\theta \mapsto q(\theta)+\phi_\gd(\theta)$ is a parametrization of $M^\gd$ (possibly $\Vert q'(\theta)+\phi_\gd'(\theta)\Vert\neq 1$).
In the following Lemma, we give the first order of the mapping $\phi_\gd$. 



\medskip

\begin{lemma}
\label{lem:first order phi}
For all $\gd$ small enough,
\begin{equation}
\label{eq: lem:first order phi}
 \sup_{\theta\in\bbR/L\bbR} \Vert\phi_\gd(\theta)-\gd h^1(\theta)\Vert\, =\, O(\gd^2)\, ,
\end{equation}
where for all $\theta\in \bbR/L\bbR$ the vector $h^1(\theta)$ is the unique solution of (recall that $H$ denotes the Hessian matrix of $V$)
\begin{equation}
 \langle h^1(\theta),q'(\theta)\rangle \, =\, 0 \qquad \text{and} \qquad H[q(\theta)]h^1(\theta)\, =\,G[q(\theta)]-\langle G[q(\theta)],q'(\theta)\rangle q'(\theta)\, .
\end{equation}
\end{lemma}

\medskip

\begin{proof}
Let $Y^\gd_0=q(\theta_0)+\phi_\gd(\theta_0)\in M^\gd$ be the initial condition of a the trajectory $Y^\gd$ of \eqref{eq:perturbed1 determ}. If we denote $\theta_t:=p(Y_t)$, then \eqref{eq:perturbed1 determ} at time $t=0$ in this case becomes
\begin{equation}
\label{eq:perturbed system on Mgd}
 (q'(\theta_0)+\phi'_\gd(\theta_0))\dot\theta^\gd_0\, =\, -\nabla V[q(\theta_0)+\phi_\gd(\theta_0)]+\gd G[q(\theta_0)+\phi_\gd(\theta_0)]\, .
\end{equation}
We view here $\dot\theta^\gd$ as a function of $\theta_0$, and we first look for uniform estimations of $\dot\theta^\gd_0$ with respect to $\theta_0$. After a projection on the tangent space of $M$ we get
\begin{equation}
\label{eq:proj for theta0}
 (1+\langle \phi'_\gd(\theta_0),q'(\theta_0)\rangle)\dot\theta^\gd_0\, =\, \langle -\nabla V[q(\theta_0)+\phi_\gd(\theta_0)]+\gd G[q(\theta_0)+\phi_\gd(\theta_0)],q'(\theta_0)\rangle\, .
\end{equation}
Recalling Lemma \ref{lem:p} and Theorem \ref{lem:phi Mgd} we deduce that $\dot \theta^\gd_0$ is $C^2$ with respect to $\theta_0$, and we get the first order expansion (using \eqref{eq:zero hess})
\begin{equation}
\dot \theta^\gd_0\, =\, \gd \langle G[q(\theta_0)],q'(\theta^\gd_0)\rangle +O(\Vert \phi_\gd(\theta_0)\Vert,\Vert\phi'_\gd(\theta_0)\Vert)\, .
\end{equation}
So we deduce from Theorem \ref{lem:phi Mgd}
\begin{equation}
\label{eq:estimate dot theta}
 \sup_{\theta_0\in\bbR/L\bbR} |\dot \theta^\gd_0|\, =\, O(\gd)\, .
\end{equation}
Now we can prove Lemma \ref{lem:first order phi}: projecting \eqref{eq:perturbed system on Mgd} on the normal space we get
\begin{multline}
 \dot\theta^\gd_0\Big(\phi'_\gd(\theta_0)-\langle \phi'_\gd(\theta_0),q'(\theta_0)\rangle q'(\theta_0)\Big)\, =\, -H[q(\theta_0)]\phi_\gd(\theta_0)+ \gd G[q(\theta_0)]\\-\gd\langle  G[q(\theta_0)],q'(\theta_0)\rangle q'(\theta_0) \\-\Big(\nabla V[q(\theta_0)+\phi_\gd(\theta_0)]-\langle \nabla V[q(\theta_0)+\phi_\gd(\theta_0)],q'(\theta_0)\rangle q'(\theta_0)-H[q(\theta_0)]\phi_\gd(\theta_0)\Big)+O(\gd^2)\, .
\end{multline}
The last line in the previous equation is of order $\gd^2$, due to Lemma \ref{lem:first order phi}, and thus for $h^1$ defined as in the statement of the Lemma we have (recall \eqref{eq:estimate dot theta})
\begin{equation}
\label{eq:first bound first order phi}
 H[q(\theta_0)](\phi_\gd(\theta_0)-\gd h^1(\theta_0))\, =\, O(\gd^2)\, .
\end{equation}
Since both vectors $\phi_\gd(\theta_0)$ and $h^1(\theta_0)$ belong to the normal space of $M$ at the point $q(\theta_0)$, the spectral gap \eqref{eq:spect gap hess} together with \eqref{eq:first bound first order phi} imply
\begin{equation}
 \phi_\gd(\theta_0)-\gd h^1(\theta_0)\, =\, O(\gd^2)\, .
\end{equation}
By a compactness argument the $O(\gd^2)$ in the previous equation is in fact uniform with respect to $\theta$, so we get \eqref{eq: lem:first order phi}.
\end{proof}

\end{subsection}
 
\begin{subsection}{Choice of projection}\label{sec:choice proj}

The proof of Theorem \ref{th:exp w} we develop is based on perturbation arguments around the manifold $M$. We will thus use the orthogonal projection on $M$ rather than the one on $M^\gd$: for a
point $Y$ located in a neighborhood of $M$, we will use the coordinates $(\theta,h)$ defined as follows
\begin{align}
 \theta\, :=\, p(Y)\, ,
\end{align}
\begin{equation}
  h\, :=\, Y-q(\theta)\, .
\end{equation}
We will use the notations $\theta^\gd_t$ and $h^\gd_t$ for a path $Y^\gd_t$ depending on $\gd$. We stress that these coordinates satisfy
\begin{equation}
 \langle h,q'(\theta)\rangle\, =\, 0\, .
\end{equation}

\begin{figure}[hlt]
\begin{center}
\leavevmode
\epsfxsize =10 cm
\psfragscanon
\psfrag{Mgd}{$M^\gd$}
\psfrag{M}{$M$}
\psfrag{Yt}{$Y$}
\psfrag{q1}{$q(\theta)$}
\psfrag{q2}{$\tilde q_\gd(\theta)=q_\gd(\tilde p_\gd(Y))$}
\psfrag{q3}{$q_\gd(\varphi)$}
\epsfbox{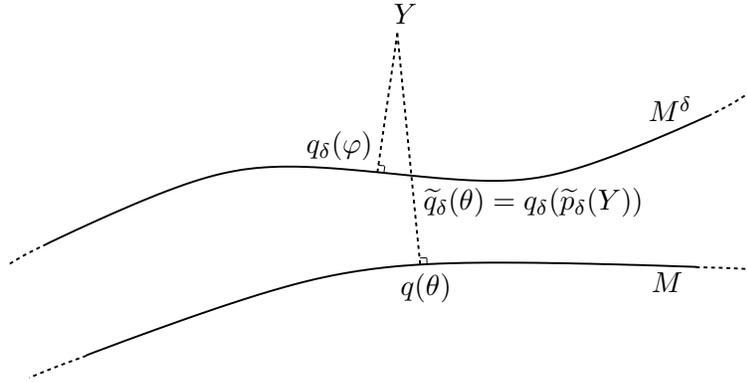}
\end{center}
\caption{$M^\gd$ parametrized by $\tilde q_\gd$. Here $\varphi=p_\gd(Y)$ and $\theta=p(Y)$.}
\label{fig:proj}
\end{figure}

We define
\begin{equation}
\label{eq:def tilde q}
 \tilde q_\gd(\theta)=q(\theta)+\phi_\gd(\theta)\, .
\end{equation} 
The parametrization $\{\tilde q_\gd(\theta),\, \theta\in \bbR/L\bbR\}$ is close in a certain sense to the one given by $q_\gd$. In fact if we define, for each point $Y$ in the neighborhood of $M^\gd$, $\tilde p_\gd(Y)$ as the phase $\varphi\in \bbR/L_\gd\bbR$ satisfying $q_\gd(\varphi)=\tilde q_\gd(\theta)$, where $\theta=p(Y)$, then we have the following Lemma:
\begin{lemma}
\label{lem:p close to tilde p}
For $\gd$ small enough and $Y$ in a neighborhood of $M^\gd$ (such that $p$ and $p_\gd$ are well defined)
\begin{equation}
\vert \tilde p_\gd(Y)- p_\gd(Y)\vert \, =\, O(\gd  \dist(Y,M^\gd))\, . 
\end{equation}
\end{lemma}

\medskip

\begin{proof}
We denote $\varphi:=p_\gd(Y)$ (recall $\theta=p(Y)$), and $\ga:=\tilde p_\gd(Y)-\varphi$. We have thus $q_\gd(\varphi+\ga)= \tilde q_\gd(\theta)$, and Theorem \ref{lem:phi Mgd} implies $q'_\gd(\varphi+\ga)=q'(\theta)+O(\gd)$. From the identity $\langle Y-q_\gd(\varphi+\ga),q'(\theta)\rangle=0$ we thus get $\langle Y-q_\gd(\varphi+\ga),q_\gd'(\varphi+\ga)+O(\gd)\rangle=0$. Expanding this last equation with respect to $\ga$ we obtain
\begin{equation}
 \langle Y-q_\gd(\varphi)-\ga q_\gd'(\varphi)+O(\ga^2),q'_\gd(\varphi)+\ga q''_\gd(\varphi)+O(\ga^2)+O(\gd)\rangle\, =\, 0\, .
\end{equation}
Expanding the scalar product, recalling the identities $\langle Y-q_\gd(\varphi),q'_\gd(\varphi)\rangle=0$ and $\Vert q'_\gd(\varphi)\Vert=1$ this reduces to
\begin{equation}
-\ga(1-\langle Y-q_\gd(\varphi),q''_\gd(\varphi)\rangle)+O(\ga^2)+O(\Vert Y-q_\gd(\varphi)\Vert \gd)\, =\, 0\, , 
\end{equation}
which implies the expected bound for $\ga$.
\end{proof}

\medskip

It will be useful to consider the restriction of the dynamics \eqref{eq:perturbed1 determ} on $M^\gd$ with respect to the parametrization
$\theta\mapsto \tilde q_\gd(\theta)$, and thus we introduce the function $\tilde b_\gd$ defined on $\bbR / L\bbR$ by
\begin{equation}
 \tilde b_\gd(\theta)\, =\, \bigg\langle -\nabla V[\tilde q_\gd(\theta)]+G[\tilde q_\gd(\theta)], \frac{{\tilde q_\gd}^\prime(\theta)}{\Vert \tilde q'_\gd(\theta)\Vert}\bigg\rangle \, .
\end{equation}
This drift satisfies the following Lemma:
\begin{lemma}
\label{lem:tilde b}
\begin{equation}
 \sup_{\theta\in \bbR/L\bbR} \{ |b(\theta)-\tilde b_\gd(\theta)|,\, |b'(\theta)-\tilde b'_\gd(\theta)|\}\, =\, O(\gd^2)
\end{equation}
\end{lemma}

This Lemma allows us to study the escape problem on tubes $U^\gd$ defined on intervals $[\varphi_{A^\gd}-\Delta_1,\varphi_{A^\gd}+\Delta_2]$ with constant length with respect to $\gd$: we will suppose in the rest of the paper that $b$ has a stable hyperbolic fixed point $\theta^0$ with domain of attraction $I$, and with this hypothesis the previous Lemma ensures that $\tilde b_\gd$ has a stable fixed $\theta^\gd_0$ located in a $\gd$-neighborhood of $\theta^0$ and whose domain of attraction is a $\gd$-perturbation of $I$. Since (recall Theorem \ref{lem:phi Mgd}) $\Vert \tilde q'_\gd(\theta)\Vert=1+O(\gd)$, using the parametrization $\theta\mapsto \tilde q'_\gd(\theta)$ instead of $q_\gd(\varphi)$ only induces an error of order $\gd$ in the phases. Thus if $\varphi_{A^\gd}$ denotes the phase satisfying $q_\gd(\varphi_{A^\gd})=\tilde q_\gd(\theta^\gd_0)$, then $\varphi_{A^\gd}$ is an hyperbolic fixed point for $b_\gd$ and its domain of attraction is also a $\gd$-perturbation of $I$. Of course, if we denote 
$A^\gd=q_\gd(\varphi_{A^\gd})$, then since $M^\gd$ is a SNHM, $A^\gd$ is a stable fixed point for \eqref{eq:perturbed1 determ}.

\medskip

\begin{proof}
 Using Theorem \ref{lem:phi Mgd} and \eqref{eq:zero grad}, it is clear that $|b(\theta)-\tilde b_\gd(\theta)|$ is of order $\gd^2$. Taking the derivative with respect to $\theta$, we obtain
\begin{multline}
\label{eq:expand bgd minus b}
 \tilde b_\gd'(\theta)-b'(\theta)\, =\, \bigg\langle -H[\tilde q_\gd(\theta)]\tilde q'_\gd(\theta),\frac{\tilde q'_\gd(\theta)}{\Vert \tilde q'_\gd(\theta)\Vert }\bigg\rangle
 +\bigg\langle -\nabla V[\tilde q_\gd(\theta)],\frac{\tilde q''_\gd(\theta)}{\Vert \tilde q'_\gd(\theta)\Vert }-\frac{\langle \tilde q''_\gd(\theta),\tilde q'_\gd(\theta)\rangle\tilde q'_\gd(\theta)}{\Vert \tilde q'_\gd(\theta)\Vert^3}\bigg\rangle
 \\ +\gd\bigg(\bigg\langle DG[\tilde q_\gd(\theta)]\tilde q'_\gd(\theta),\frac{\tilde q'_\gd(\theta)}{\Vert \tilde q'_\gd(\theta)\Vert }\bigg\rangle-\langle DG[q(\theta)]q'(\theta),q'(\theta)\rangle \bigg)
 \\ +\gd \bigg(\bigg\langle G[\tilde q_\gd(\theta)],\frac{\tilde q''_\gd(\theta)}{\Vert \tilde q'_\gd(\theta)\Vert }-\frac{\langle \tilde q''_\gd(\theta),\tilde q'_\gd(\theta)\rangle\tilde q'_\gd(\theta)}{\Vert \tilde q'_\gd(\theta)\Vert^3}\bigg\rangle-\langle G[q(\theta)],q''(\theta)\rangle\bigg)\, .
\end{multline}
Using Theorem \ref{lem:phi Mgd} we get the following expansion for the first term of the right hand side (recall that it implies in particular $\Vert \tilde q'_\gd(\theta)\Vert=1+O(\gd)$):
\begin{multline}
 \bigg\langle H[\tilde q_\gd(\theta)]\tilde q'_\gd(\theta),\frac{\tilde q'_\gd(\theta)}{\Vert \tilde q'_\gd(\theta)\Vert }\bigg\rangle\, =\, \langle H[q(\theta)]q'(\theta),q'(\theta)\rangle+\langle H[q(\theta)](\tilde q_\gd'(\theta)-q'(\theta)),q'(\theta)\rangle \\
 +\bigg\langle H[q(\theta)] q'(\theta), \frac{\tilde q'_\gd(\theta)}{\Vert \tilde q'_\gd(\theta)\Vert }-q'(\theta)\bigg\rangle+D^3V[q(\theta)](\phi_\gd(\theta),q'(\theta),q'(\theta))+O(\gd^2)\, ,
\end{multline}
and \eqref{eq:zero hess} implies that the three first term of the right hand side in this expansion are equal to zero. Using similar argument to treat the other terms 
of \eqref{eq:expand bgd minus b} (recalling in particular \eqref{eq:zero grad}), we see that it reduces to
\begin{equation}
 \tilde b_\gd'(\theta)-b'(\theta)\, =\,-D^3 V[q(\theta)](\phi_\gd(\theta),q'(\theta),q'(\theta))-\langle H[q(\theta)]\phi_\gd(\theta),q''(\theta)\rangle+O(\gd^2)\, .
\end{equation}
Now two derivations with respect to $\theta$ of the identity $\nabla V[q(\theta)]=0$ imply that for all $u\in \bbR^n$
\begin{equation}
\label{eq:link V3 H}
 D^3V[q(\theta)](q'(\theta),q'(\theta),u)+\Big\langle H[q(\theta)] q''(\theta),u\Big\rangle\, =\, 0\, ,
\end{equation}
and for $u=q'(\theta)$ this implies the expected bound for $\tilde b_\gd'(\theta)-b'(\theta)$.
\end{proof}

\end{subsection}

\end{section}

\medskip

\begin{section}{Quasipotential and optimal path}\label{sec:optimal path}

 As shown in \cite{cf:DayDarden}, a continuity argument allows us to define $W_\gd(A^\gd,\cdot)$ as the infimum of the rate function over the paths defined on $(-\infty,0]$ with limit $A^\gd$ at $-\infty$. In fact, extending the paths $Y\in C([-T,0],\bbR^n)$ with $Y_{-T}=A^\gd$ by $Y_t=A^\gd$ for $t\leq -T$, we get that for all $E\in \bbR^n$
\begin{equation}
\label{ineq:W infinite time}
 W_\gd(A^\gd,E)\, \geq\, \inf\left\{I_{\gd,-\infty}^{A^\gd}(Y):\,Y\in C((-\infty,0],\bbR^n),\, \lim_{t\rightarrow -\infty} Y_t=A^\gd,\, Y_{0}=E\right\}\, .
\end{equation}
On the other hand, for each path $Y\in C((-\infty,0],\bbR^n)$ with  $\lim_{t\rightarrow -\infty} Y_t=A^\gd$ and $Y_{0}=E$, we have for all $t\leq 0$
\begin{equation}
 W_\gd(A^\gd,E)\, \leq \, W_\gd(A^\gd,Y_{t})+I_{\gd,t}^{Y_{t}}(Y)\, .
\end{equation}
But $W_\gd(A,\cdot)$ is Lipschitz continuous (see \cite{cf:FreidWentz} Lemma 2.3), so $W_\gd(A^\gd,Y_t)\rightarrow 0$ when $t \rightarrow -\infty$, and thus \eqref{ineq:W infinite time} is in fact an equality. 

\medskip

For a point $E\in \bbR^n$, we call an optimal path a path $Y^\gd\in C((-\infty,0],\bbR^n)$ with $\lim_{t\rightarrow -\infty} Y^\gd_t= A^\gd$, $Y_0=E$ and satisfying
\begin{equation}
 I_{\gd,-\infty}^{A^\gd}(Y^\gd)\, =\, W_\gd(A^\gd,E)\, .
\end{equation}
In \cite{cf:DayDarden} it is explained that for each $E\in \bbR^n$, if the trajectories approximating $W_\gd(A^\gd,E)$ stay in a compact,
then a compactness argument ensures the existence of an optimal path for $E$. We follow this idea in two steps, using in the second step the hyperbolic structure of $M$.

\begin{lemma}
\label{lem:opt path min}
Let $K$ be a compact neighborhood of $A^\gd$.
Then there exists an optimal path for at least one point $B^\gd\in \partial K$ satisfying
\begin{equation}
\label{eq:min Wgd Bgd K}
 W_\gd(A^\gd,B^\gd)\, =\, \inf_{E\in \partial K}W_\gd(A^\gd,E)\, .
\end{equation}
\end{lemma}

\medskip

\begin{proof}
We can choose a sequence of paths $Y^{k,\gd}$ staying in $K$ with $Y^{k,\gd}_0\in \partial K$, such that $I_{\gd,-\infty}^{A^\gd}(Y^{k,\gd})$ converges to $\inf_{E\in \partial K} W_\gd(A^\gd,E)$. Using the Arzel\`a-Ascoli theorem (the convergence of the rate function and the compactness ensures a uniform control of $\int_t^{t+\gep}\Vert \dot Y^{k,\gd}_s \Vert^2\dd s$ and thus the equicontinuity) on each compact interval of time $[-m,0]$ and a diagonal procedure we can show that there exists a sub-sequence $Y^{\psi_k,\gd}$ that converges in $C((-\infty,0],\bbR^n)$ to a path $Y^\gd$, which has the expected properties.
\end{proof}

\medskip

The previous Lemma does not give the existence of an optimal path for each $B^\gd \in \partial U^\gd$ (recall \eqref{eq:def tube U}) satisfying \eqref{eq:min Wgd Bgd th}. To get this result, we rely on the hyperbolic
structure of $M$ to ensure that the trajectories approximating $W_\gd(A^\gd,B^\gd)$ stay in a compact, for each of these points. 

\medskip

\begin{lemma}
\label{lem:opt path}
For $\gd$ small enough there exists an optimal path $Y^\gd$ for each $B^\gd \in \partial U^\gd$ satisfying \eqref{eq:min Wgd Bgd th}. Moreover if we define the whole tube
\begin{equation}
 \tilde U^\gd\, =\, \{Z\in \bbR^n:\, \dist(Z,M^\gd)\leq C_0 \gd^{1/2}\}\, ,
\end{equation}
then $Y^\gd_t \in \tilde U^\gd$ for all $t\leq 0$.
\end{lemma}

\medskip

\begin{proof}

We first aim at proving that a trajectory of \eqref{eq:perturbed1 determ} starting on the boundary $\partial\tilde U^\gd$ is strictly inside $\tilde U^\gd$ for times small enough. To show this property, we just need to prove that for all $Z\in \partial \tilde U^\gd$
\begin{equation}
\label{ineq:traj bound tilde U}
 \langle -\nabla V[Z]+\gd G[Z], Z-q_\gd(p_\gd(Z))\rangle \, <\, 0\, .
\end{equation}
Now if we denote $\theta=p(Z)$, Lemma \ref{lem:p close to tilde p} implies (recall \eqref{eq:def tilde q} and $dist(Z,M^\gd)= C_0 \gd^{1/2}$)
\begin{equation}
 q_\gd(p_\gd(Z))\, =\, \tilde q_\gd(\theta)+O(\gd^{3/2}),
\end{equation}
but in fact a control of order $\gd^2$ for the error term is enough for our purpose, and recalling Theorem \ref{lem:phi Mgd} we obtain
\begin{equation}
 q_\gd(p_\gd(Z))\, =\, q(\theta)+O(\gd)\, .
\end{equation}
This implies in particular
\begin{equation}
\label{eq:norm Z-q}
 \Vert Z-q(\theta)\Vert \, =\, C_0\gd^{1/2}+O(\gd)\, .
\end{equation}
Now recalling \eqref{eq:zero grad} we get the first order expansion
\begin{equation}
 \nabla V[Z]\, =\, H[q(\theta)](Z-q(\theta))+O(\gd)\, .
\end{equation}
We deduce
\begin{equation}
  \langle -\nabla V[Z]+\gd G[Z], Z-q_\gd(p_\gd(Z))\rangle\, =\, -\langle H[q(\theta)](Z-q(\theta)),Z-q(\theta)\rangle +O(\gd^{3/2})\, ,
\end{equation}
and the spectral gap \eqref{eq:spect gap hess} implies that the first term of the right hand side is bounded from above by $-\gl \Vert Z-q(\theta)\Vert^2$, so (recall \eqref{eq:norm Z-q}) \eqref{ineq:traj bound tilde U} is satisfied for $\gd$ small enough. 
Now define the compact
\begin{equation}
 K\, :=\, \{Z\in\bbR^n:\, \dist(Z,M^\gd)\leq (C_0+1)\gd^{1/2}\}\, . 
\end{equation}
From Lemma \ref{lem:opt path min}, we know that there exists an optimal path $Y^\gd$ for a point $Q^\gd\in \partial K$ satisfying $W_\gd(A^\gd,Q^\gd)\, =\, \inf_{E\in \partial K}W_\gd(A^\gd,E)$. Define
 \begin{equation}
  T\, =\, \sup\{t: \, Y^\gd_{t}\in \tilde U^\gd\}\, .
 \end{equation}
Then $T<0$ and we have
\begin{equation}
 \inf_{E\in \partial K} W_\gd(A^\gd,E)\, \geq \, \inf_{E\in \partial \tilde U^\gd} W_\gd(A^\gd, E)+ I_{\gd,T}^{Y^\gd_{T}}(Y^\gd)\, . 
\end{equation}
\eqref{ineq:traj bound tilde U} shows that the last term of the previous equation is positive, so $\inf_{E\in \partial \tilde U^\gd} W_\gd(A^\gd, E)<\inf_{E\in \partial K} W_\gd(A^\gd,E)$. Consequently, for a point $B^\gd \in \partial U^\gd$ satisfying \eqref{eq:min Wgd Bgd th}, a trajectory $Z^\gd$ linking $A^\gd$ to $B^\gd$ in such a way that $I_{\gd,-\infty}^{A^\gd}(Z^\gd)$ is sufficiently close to $W_\gd(A^\gd,B^\gd)$ must stay in $K$. So an argument similar to the proof of Lemma \ref{lem:opt path min}, involving compactness, ensures the existence of an optimal path for $B^\gd$.

If such an optimal path, say $Y^\gd$, exits $\tilde U^\gd$ at a time $t_0<0$, then $I_{\gd,-\infty}^{A^\gd}(Y^\gd)\geq W_\gd(A^\gd,Y^\gd_{t_0})+I_{\gd,t_0}^{Y^\gd_{t_0}}(Y^\gd)$. Since $Y^\gd$ is optimal and $B^\gd$ satisfies \eqref{eq:min Wgd Bgd th}, this implies $I_{\gd,t_0}^{Y^\gd_{t_0}}(Y^\gd)=0$, i.e the remaining part of the path $Y^\gd$ is solution of \eqref{eq:perturbed1 determ}, and it contradicts \eqref{ineq:traj bound tilde U}.

\end{proof}

\medskip

\begin{rem}\label{rem:opt path tilde U}
 The preceding proof can be easily adapted to show that there exists an optimal path staying in $\tilde U^\gd$ for each point $B^\gd\in \partial\tilde U^\gd$ satisfying
 \begin{equation}
  W_\gd(A^\gd,B^\gd)\, =\, \inf_{E\in \partial\tilde U^\gd}W_\gd(A^\gd,E)\, .
 \end{equation}
This will be useful to prove the existence of an optimal path for each point of the type $q_\gd(\varphi^\gd)$ in the proof of Corollary \ref{cor:w} in Subsection \ref{sec:proof cor}.
\end{rem}

\medskip

Lemma \ref{lem:opt path} ensures that the optimal paths of each point $B^\gd\in\partial U^\gd$ satisfying \eqref{eq:min Wgd Bgd th} stay in the whole tube $\tilde U^\gd$, but it does not ensure that they stay in the truncated one $U^\gd$. Unfortunately we are not able to prove directly this fact. It will be a consequence of the following preliminary Lemma and Lemma \ref{lem:exit} (see remark \ref{rem:dot theta h 1/2}). But it does not cause any problem for the proofs preceding Lemma \ref{lem:exit}, since they only lean on the fact that the optimal path stay in a $\gd^{1/2}$-neighborhood of $M^\gd$ (which is of course induced by Lemma \ref{lem:opt path}).

\medskip

\begin{lemma}
\label{lem:prelim exit Y}
For all $\gep>0$ there exists $\gd_\gep>0$ such that if $\gd\leq\gd_\gep$, for all $B^\gd \in \partial U^\gd$ satisfying \eqref{eq:min Wgd Bgd th} and all associated optimal path $Y^\gd$, if $Y^\gd$ leaves $U^\gd$ at a time $t_0<0$, then
\begin{equation}
 \dist(Y^\gd_{t_0},M^\gd)\, \geq\, \gd^{1/2+\gep}\, .
\end{equation}
\end{lemma}

\medskip

\begin{proof}

Consider a $Z\in \partial U^\gd$ such that $\Vert Z-q_\gd(p_\gd(Z))\Vert \leq \gd^{1/2+\gep}$. We have for $\theta=p(Z)$, proceeding as in the previous Lemma (recalling Lemma \ref{lem:p close to tilde p}),
\begin{multline}
 \langle -\nabla V[Z]+\gd G[Z],q_\gd'(p_\gd(Z))\rangle \\ \, =\, \Big\langle -\nabla V[q(\theta)]-H[q(\theta)](Z-q(\theta))+\gd G[q(\theta)]+O(\gd^{1+2\gep},\gd^{3/2+\gep}),q'(\theta)+\phi'_\gd(\theta)+O(\gd^{3/2+\gep})\Big\rangle\, .
\end{multline}
So using \eqref{eq:zero grad}, \eqref{eq:zero hess}, and Theorem \ref{lem:phi Mgd} we obtain (recalling \eqref{eq:def b})
\begin{equation}
 \langle -\nabla V[Z]+\gd G[Z],q_\gd'(p_\gd(Z))\rangle \, =\, b(\theta)+O(\gd^{1+2\gep},\gd^{3/2+\gep})\, .
\end{equation}
We deduce (reminding that $[\varphi_{A^\gd}-\Delta_1,\varphi_{A^\gd}+\Delta_2]$ is included in the domain of attraction of $\theta^0$ for $b$) that the trajectories of \eqref{eq:perturbed1 determ} starting at such points $Z$ are strictly in $U^\gd$ for times small enough. A premature exit of an optimal path $Y^\gd$ can not occur at such a point.

\end{proof}

\medskip

We now give the Euler Lagrange type equation satisfied by the optimal paths. It corresponds to Theorem 1 in \cite{cf:DayDarden}. We denote by $A^\dagger$ the transpose of a square matrix $A$.

\begin{lemma}
\label{lem:EL}
Let $E\in\bbR^n$ admitting an optimal path $Y^\gd$. Then $Y^\gd\in C^2((-\infty,0),\bbR^n)$ and satisfies for all $t<0$
\begin{equation}
 \label{eq:EL}
 \ddot Y^\gd_t\, =\, \left(H[Y^\gd_t]-\gd  DG^\dagger[Y^\gd_t]\right)\left(\nabla V[Y^\gd_t]-\gd G[Y^\gd_t]\right)+\gd \left(DG[Y^\gd_t]-DG^\dagger[Y^\gd_t]\right)\dot Y^\gd_t\, .
\end{equation}
\end{lemma}

\medskip

\begin{proof}
Define
\begin{equation}
 I_{\gd,T_1,T_2}(Y)\, =\, \frac{1}{2}\int_{T_1}^{T_2}\Vert\dot Y_t+ \nabla V[Y_t]-\gd G[Y_t]  \Vert^2\dd t \, .
\end{equation}
For all $T_1<T_2<0$ an optimal path $Y^ \gd$ must be a local minimum for $I_{\gd,T_1,T_2}(Z)$ viewed as a functional on the space of absolutely continuous paths $Z$ satisfying $Z_{T_1}=Y^\gd_{T_1}$ and $Z_{T_2}=Y^\gd_{T_2}$. We denote $H^k((T_1,T_2),\bbR^n)$ the usual Sobolev spaces on the interval $(T_1,T_2)$. Remark that since $I_{\gd,T_1,T_2}(Y^\gd)<\infty$, $Y^\gd_{|(T_1,T_2)}\in H^1((T_1,T_2),\bbR^n)$, and in particular the right hand side of \eqref{eq:EL} is well defined in the sense of distributions. Let $f\in C^\infty((T_1,T_2),\bbR^n)$ with compact support. We get the expansion for $\eta\in \bbR$
\begin{multline}
 I_{\gd,T_1,T_2}(Y^\gd+\eta f)\, =\, I_{\gd,T_1,T_2}(Y^\gd)+\eta\int_{T_1}^{T_2}\left\langle \dot Y^\gd_t,\dot f_t\right\rangle+\left\langle\nabla V[Y^\gd_t]-\gd G[Y^\gd_t],\dot f_t\right\rangle \\ 
 +\left\langle  H[Y^\gd_t] f_t -\gd DG[Y^\gd_t]f_t,\dot Y^\gd_t\right\rangle 
 +\left\langle H[Y^\gd_t]f_t-\gd DG[Y^\gd_t]f_t,\nabla V[Y^\gd_t]-\gd G[Y^\gd_t]\right\rangle \dd t \\ +O(\eta^2)\, .
\end{multline}
Since $Y^\gd$ is a local minimum, the term of order $\eta$ in the right hand side of previous equation is equal to $0$, and it implies \eqref{eq:EL} on the interval $(T_1,T_2)$ in the sense of distributions. But since $Y^\gd_{|(T_1,T_2)}\in H^1((T_1,T_2),\bbR^n)$, $V$ is $C^4$ and $G$ is $C^2$, \eqref{eq:EL} implies that $\ddot Y^\gd_{|(T_1,T_2)}\in L^2((T_1,T_2),\bbR^n)$, or in other words $Y^\gd_{|(T_1,T_2)}\in H^2((T_1,T_2),\bbR^n)$. But again in this case \eqref{eq:EL} implies that $\ddot Y^\gd_{|(T_1,T_2)}\in H^1((T_1,T_2),\bbR^n)$, and thus admits a continuous representation.
\end{proof}

\end{section}

\medskip

\begin{section}{Proof of Theorem \lowercase{\ref{th:exp w}} and Corollary \lowercase{\ref{cor:w}}}

\begin{subsection}{Sketch of the proof}

The aim of the proof is to make a expansion of $I_{-\infty,\gd}^{A^\gd}(Y^\gd_t)$ for the optimal paths linking $A^\gd$ to the points $B^\gd\in \partial U^\gd$ satisfying 
\eqref{eq:min Wgd Bgd th}, and to compare this expansion to \eqref{eq:I reduced}. The main idea we follow is that when a trajectory $Y_t$ is located at a distance of order $\gd^2$ from $M^\gd$ on a time interval $[T_1,T_2]$, then it is possible to make an accurate expansion of 
\begin{equation}
\label{I T1 T2}
 \int_{T_1}^{T_2}\Vert \dot Y_t+\nabla V[Y_t]-\gd G[Y_t]\Vert^2\dd t\, .
\end{equation}
Unfortunately we are not able to prove that a optimal trajectory linking $A^\gd$ to a point $B^\gd\in \partial U^\gd$ satisfying \eqref{eq:min Wgd Bgd th} stays at distance $\gd^2$ from $M^\gd$. However we are able to prove that such a optimal path stays at distance of order $\gd^2$ for most of the time. To do that we rely on the fact that the optimal paths satisfy the Euler Lagrange type equation \eqref{eq:EL}. When $\gd=0$ \eqref{eq:EL} reduces to
\begin{equation}
\label{eq:EL gd null}
 \ddot Y_t\,=\, H[Y_t]\nabla V[Y_t]\, ,
\end{equation}
and a solution of \eqref{eq:EL gd null} starting in a neighborhood of $M$ but not in $M$ moves away from $M$. Indeed if $q(\theta)$ is the projection on $M$ of a point $Y$ located in a neighborhood of $M$, then
\begin{equation}
 H[Y]\nabla V[Y]\, =\, H[Y]^2 (Y-q(\theta))+O(\Vert Y-q(\theta)\Vert^2)\, ,
\end{equation}
and \eqref{eq:zero hess} and the spectral gap \eqref{eq:spect gap hess} imply that the vector $H[Y]^2 (Y-q(\theta))$ lies in the normal space of $M$ at the point $q(\theta)$ and has a norm bounded from below by $\gl^2 \Vert Y-q(\theta)\Vert$.
Using perturbation arguments, we show that when $\gd\neq 0$, a solution of the Euler Lagrange type equation \eqref{eq:EL} starting from a point located at a distance from $M^\gd$ bounded from below by $C\gd^2$ behaves similarly: it moves away from $M^\gd$. So for all optimal path $Y^\gd$ there exists a time $\tau^\gd_1\leq 0$ (maybe equal to $0$) such that for $t\leq \tau^\gd_1$ $Y^\gd$ is located at distance of order $\gd^2$ from $M^\gd$, and if $\tau^\gd_1<0$ then for $t>\tau^\gd_1$ $\dist (Y^\gd_t,M^\gd)\geq C\gd^2$. It is the purpose of Lemma \ref{lem:dist gd2 M}, where we also control $|\tau^\gd_1|$ when it is non null. The Lemmas \ref{lem:dot Y 1/2}, \ref{lem:dist gd M} and \ref{lem:dot Y gd} are intermediate results needed to prove Lemma \ref{lem:dist gd2 M}. In Lemma \ref{lem:dot h 2} we control the derivate in time of $Y^\gd$ for $t\leq \tau^\gd_1$.

In Lemma \ref{lem:exit} we control the behavior of the optimal paths $Y^\gd$ on the time interval $[\tau^\gd_1,0]$ (if $\tau^\gd_1 <0$). The estimations we obtain are sufficient to allow a good expansion of
\begin{equation}
 \int_{\tau^\gd_1}^{0}\Vert \dot Y^\gd_t+\nabla V[Y^\gd_t]-\gd G[Y^\gd_t]\Vert^2\dd t\, .
\end{equation}

Finally the expansion of \eqref{I T1 T2} we are able to make depends on the length of the time interval $[T_1,T_2]$, so we can not simply take $T_1=-\infty$ and $T_2=\tau^\gd_1$. We have to restrict the expansion on a time interval $[\tau^\gd_0,\tau^\gd_1]$, choosing $\tau^\gd_0$ in such a way that
$\int_{-\infty}^{\tau^\gd_0}\Vert \dot Y^\gd_t+\nabla V[Y^\gd_t]-\gd G[Y^\gd_t]\Vert^2\dd t$ is negligible.

\begin{figure}[hlt]
\begin{center}
\leavevmode
\epsfxsize =13.5 cm
\psfragscanon
\psfrag{M}{$M^\gd$}
\psfrag{B}{$B^\gd$}
\psfrag{A}{$A^\gd$}
\psfrag{l}{$\sim \gd^2$}
\psfrag{L}{$\sim \gd|\log\gd|$}
\psfrag{X}{$Y^\gd_{\tau^\gd_0}$}
\psfrag{Y}{$Y^\gd_{\tau^\gd_1}$}
\psfrag{U}{$U^\gd$}
\psfrag{lu}{$\sim \gd^{1/2}$}
\epsfbox{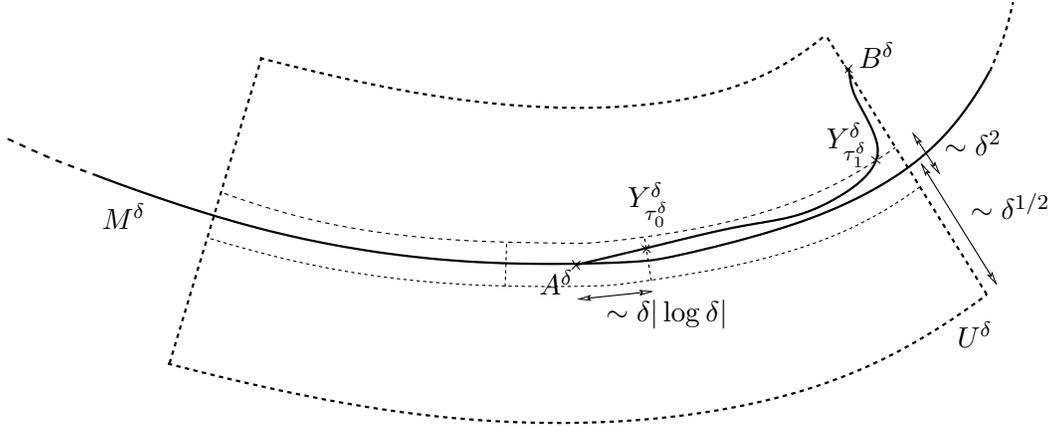}
\end{center}
\caption{The times $\tau^\gd_0$ and $\tau^\gd_1$ for an optimal path $Y^\gd$ linking $A^\gd$ to $B^\gd$.}
\label{fig:heuristique}
\end{figure}

\end{subsection}

\begin{subsection}{Preliminary results}




We can easily find a first upper bound for $W^{red}_\gd(\varphi_{A^\gd},\varphi)$ for all $\varphi$: we deduce indeed from Theorem \ref{lem:phi Mgd} that there exists $C>0$ such that for $\gd$ small enough $|b_\gd|\leq C\gd$ (recall the definition of $\gd$ \eqref{eq:def b gd}) and thus, since the lengths of the curves $M^\gd$ are also bounded, there exists $C_1>0$ such that for all $\gd$ small enough
\begin{equation}
\label{ineq:first bound Wred}
 \sup_{\varphi\in \bbR/L_\gd \bbR}W^{red}_\gd(\varphi_{A^\gd},\varphi) \, \leq\, C_1 \gd\, .
\end{equation}
Recalling \eqref{ineq:comparaison W Wred} it implies in particular that for $B^\gd\in \partial U^\gd$ satisfying \eqref{eq:min Wgd Bgd th} we have
\begin{equation}
 \label{ineq:first bound W}
 W_\gd(A^\gd,B^\gd)\,\leq \, C_1\gd\, .
\end{equation}

\medskip

The spectral gap \eqref{eq:spect gap hess} implies that for points $E$ sufficiently close to $M$ we have
\begin{equation}
 V[E]\, \geq\, \frac{\gl}{2}\dist(E,M)\, .
\end{equation}
So we can choose the value of $C_0$ (recall the definition of $U^\gd$ in Theorem \ref{th:exp w}) such that for $\gd$ small enough
\begin{equation}
\label{eq:choice C0}
 \inf_{\dist(E,M^\gd)=C_0\gd^{1/2}}V[E]\, \geq\, 2 C_1\gd\, .
\end{equation}
This choice of $C_0$ will be useful in the proof of Lemma \ref{lem:exit} below, to ensure that the projection $p_\gd(B^\gd)$ coincides either with $\varphi_{A^\gd}-\Delta_1$ or with $\varphi_{A^\gd}+\Delta_2$.

\medskip

From Lemma \ref{lem:opt path} we know that there exists at least an optimal path $Y^\gd$ for $B^\gd$, and that $\sup_{t\leq 0}\dist(Y^\gd,M^\gd)\leq C_0\gd^{1/2}$. Since $M^\gd$ is located at a distance of order $\gd$ from $M$ (see Theorem \ref{lem:phi Mgd}), this implies (recall the coordinates introduced in Section \ref{sec:choice proj}) that $h^\gd_t$ is of order $\gd^{1/2}$ for all $t$: there exists $C>0$ such that
\begin{equation}
\label{ineq:first bound h}
 \sup_{t\leq 0}\Vert h^\gd_t\Vert\, \leq\, C\gd^{1/2}\, . 
\end{equation}
$Y^\gd$ converges to $A^\gd$ when $t\rightarrow -\infty$, so Lemma \ref{lem:EL} ensures that its second derivate is bounded on $(-\infty,0)$, and thus $\Vert \dot Y^\gd\Vert$ is also bounded and since $Y^\gd$ stays in a compact (see Lemma \ref{lem:opt path}), $\Vert \dot Y^\gd\Vert$ reaches its maximum on $(-\infty,0)$.
The following Lemmas give some properties satisfied by $Y^\gd$ and its derivative in time. We will drop the dependence in the initial value in the large deviation rate for simplicity:
\begin{equation}
 I_{\gd,T}(Z)\, :=\, I_{\gd,T}^{Z_{T}}(Z)\, .
\end{equation}
The constant $C$ is a generic constant independent from $\gd$, and whose value may change during the proof.

\medskip

\begin{lemma}
\label{lem:dot Y 1/2}
There exists $C_2>0$ such that for all $\gd$ small enough
\begin{equation}
 \sup_{t\leq 0}\Vert \dot Y^\gd_t\Vert \, \leq\, C_2 \gd^{1/2}\, .
\end{equation}
\end{lemma}

\medskip

\begin{rem}
\label{rem:dot theta h 1/2}
 Due to the smoothness of the projection $p$ (see Lemma \ref{lem:p}), the result of the previous Lemma is also true for the coordinates $\theta^\gd$ and $h^\gd$ associated with $Y^\gd$:
 \begin{equation}
  \sup_{t\leq 0} \{|\dot \theta^\gd_t|,\Vert \dot h^\gd_t \Vert\}\, =\, O(\gd^{1/2})\, .
 \end{equation}
The same argument will be true for Lemma \ref{lem:dot Y gd} below.
\end{rem}

\medskip

\begin{proof}
From Lemma \ref{lem:EL}, \eqref{eq:zero grad} and \eqref{ineq:first bound h} we deduce that there exists $C>0$ such that
\begin{equation}
\label{bound:ddot Y}
 \Vert \ddot Y^\gd_t\Vert \, \leq\,  C(\gd^{1/2}+\gd \Vert \dot Y^\gd_t\Vert)\, .
\end{equation}
Now suppose that $\sup_{t\leq 0} \Vert \dot Y^\gd_t\Vert$ is reached at $\Vert\dot Y^\gd_{t_0}\Vert$. Then, for $s\in  [t_0-1,t_0]$, $\Vert \dot Y^\gd_s\Vert$ satisfies
\begin{equation}
 \Vert\dot Y^\gd_s\Vert\, \geq\, \Vert\dot Y^\gd_{t_0}\Vert-C(\gd^{1/2}+\gd\Vert \dot Y^\gd_{t_0}\Vert)\, .
\end{equation}
So for $\gd$ small enough
\begin{equation}
\label{ineq:bound dot Y below}
 \Vert\dot Y^\gd_s\Vert\, \geq\, \frac12\Vert \dot Y^\gd_{t_0}\Vert-C\gd^{1/2}\, .
\end{equation}
Using the elementary bound $(a+b)^2\geq \frac12 a^2-b^2$, we can bound $I_{\gd,-\infty}^{A^\gd}(Y^\gd)$ from below (recall \eqref{eq:I}):
\begin{equation}
\label{bound:I below}
 I_{\gd,-\infty}(\varphi^\gd)\, \geq\, \frac14 \int_{t_0-1}^{t_0}\Vert \dot Y^\gd_t\Vert^2\dd t-\frac12\int_{t_0-1}^{t_0}\Vert \nabla V[Y^\gd_t] -\gd G[Y^\gd_t]\Vert^2\dd t\, .
\end{equation}
Taking into account \eqref{ineq:bound dot Y below}, \eqref{eq:zero grad} and \eqref{ineq:first bound h} we get
\begin{equation}
\label{bound:I below 2}
I_{\gd,-\infty}(\varphi^\gd)\, \geq\,\frac14\left(\frac12 \Vert\dot Y^\gd_{t_0}\Vert-C\gd^{1/2}\right)^2-C\gd\, .
\end{equation}
So there exists $C_2>0$ such that if $\dot Y^\gd_{t_0}\geq C_2 \gd^{1/2}$, \eqref{bound:I below 2} contradicts \eqref{ineq:first bound W}.
\end{proof}

\medskip

\begin{lemma}
\label{lem:dist gd M}
There exists $C_3>0$ such that for all $\gd$ small enough, if we define
\begin{equation}
 \tau^\gd_{2}\, :=\, \inf \{t<0, \Vert h^\gd_t-\phi_\gd(\theta^\gd_t)\Vert \geq C_3\gd\}\, ,
\end{equation}
then $\Vert h^\gd_t-\phi_\gd(\theta^\gd_t)\Vert \geq C_3\gd$ for all $ \tau^\gd_{2}\leq t\leq 0$, and if $\tau^\gd_{2}<0$ then $\Vert h^\gd_t-\phi_\gd(\theta^\gd_t)\Vert$ increases strictly on $(\tau^\gd_{2},0]$.
\end{lemma}

\medskip

\begin{proof}
Using Lemma \ref{lem:first order phi}, \eqref{eq:zero grad} and \eqref{ineq:first bound h} we get 
\begin{multline}
H[Y^\gd_t](\nabla V[Y^\gd_t]-\gd G[Y^\gd_t])\, =\, H[q(\theta^\gd_t)](H[q(\theta^\gd_t)]h^\gd_t-\gd G[q(\theta^\gd_t)])+O(\gd)\\
=\, \Big(H[q(\theta^\gd_t)]\Big)^2(h^\gd_t-\phi_\gd(\theta^\gd_t))+O(\gd)\, ,
\end{multline}
thus using furthermore Lemma \ref{lem:dot Y 1/2} we get the following first order expansion of \eqref{eq:EL}:
\begin{equation}
\label{eq:EL first order}
 \ddot Y^\gd_t\, =\, \Big(H[q(\theta^\gd_t)]\Big)^2(h^\gd_t-\phi_\gd(\theta^\gd_t))+O(\gd)\, .
\end{equation}
Now define
\begin{equation}
 \label{eq:def alpha}
 \ga^\gd_t\, =\,\Vert h^\gd_t-\phi_\gd(\theta^\gd_t)\Vert^2\, .
\end{equation}
A straightforward calculation gives
\begin{equation}
\label{eq:ddot alpha}
 \ddot \alpha^\gd_t\, =\, 2\Vert \dot h^\gd_t -\dot\theta^\gd_t\phi_\gd'(\theta^\gd_t)\Vert^2+2\langle  h^\gd_t-\phi_\gd(\theta^\gd_t),\ddot h^\gd_t-\ddot\theta^\gd_t\phi_\gd'(\theta^\gd_t)-(\dot\theta^\gd_t)^2 \phi_\gd''(\theta^\gd_t)\rangle\, .
\end{equation}
Using Lemma \ref{lem:p} we express $\dot h^\gd_t$ in with respect to $\dot Y^\gd_t$ and $h^\gd_t$:
\begin{equation}
\label{eq:dot h}
 \dot h^\gd_t\, =\, \dot Y^\gd_t -\frac{1}{1-\langle h^\gd_t,q''(\theta^\gd_t)\rangle}\langle \dot Y^\gd_t,q'(\theta^\gd_t)\rangle q'(\theta^\gd_t)\, ,
\end{equation}
and after a derivation in time it leads to
\begin{multline}
\label{eq:ddot h}
 \ddot h^\gd_t\, =\, \ddot Y^\gd_t -\dot \theta^\gd_t\frac{1}{1-\langle h^\gd_t,q''(\theta^\gd_t)\rangle}q''(\theta^\gd_t)      +\bigg(\frac{\langle \dot h^\gd_t,q''(\theta^\gd_t)\rangle+\dot \theta^\gd_t\langle h^\gd_t, q'''(\theta^\gd_t)\rangle}{(1-\langle h^\gd_t,q''(\theta^\gd_t)\rangle)^2}  \\-\frac{1}{1-\langle h^\gd_t,q''(\theta^\gd_t)\rangle}(\langle \dot Y^\gd_t,q'(\theta^\gd_t)\rangle+\dot\theta^\gd_t\langle Y^\gd_t,q''(\theta^\gd_t)\rangle\bigg)q'(\theta^\gd_t)\, .
\end{multline}
So taking together \eqref{eq:EL first order}, \eqref{eq:ddot alpha}, \eqref{eq:ddot h}, Lemma \ref{lem:dot Y 1/2} and Remark \eqref{rem:dot theta h 1/2}, we get (recall $\langle h^\gd_t,q'(\theta^\gd_t)\rangle=0$)
\begin{equation}
\label{eq:ddot ga first order 2}
 \ddot \ga^\gd_t\, =\,2\Vert \dot h^\gd_t -\dot\theta^\gd_t\phi_\gd'(\theta^\gd_t)\Vert^2+ 2\Vert H[q(\theta^\gd_t)](h^\gd_t-\phi_\gd(\theta^\gd_t))\Vert^2+O(\gd\Vert h^\gd_t-\phi_\gd(\theta^\gd_t)\Vert)\, ,
\end{equation}
where the second term of the right hand side is bounded below by $\gl^2\Vert h^\gd_t-\phi_\gd(\theta^\gd_t)\Vert^2$, due to \eqref{eq:spect gap hess}. So there exists $C_3$ such that if $\Vert h^\gd_t-\phi_\gd(\theta^\gd_t)\Vert\geq C_3\gd$ and $\dot\ga^\gd_t\geq 0$ for a $t<0$, then $\ga^\gd_s$ is strictly increasing for $s>t$.
\end{proof}

\medskip

\begin{lemma}
\label{lem:dot Y gd}
There exists $C_4>0$ such that for all $\gd$ small enough 
\begin{equation}
 \sup_{t\leq \tau^\gd_{2}}\Vert \dot Y^\gd_t\Vert\, \leq\, C_4\gd\, .
\end{equation}
\end{lemma}

\medskip

\begin{proof}
From Lemma \ref{lem:EL}, Lemma \ref{lem:dist gd M} and \eqref{eq:zero grad}  we deduce that there exists $C>0$ such that for $t\leq \tau^\gd_{2}$
\begin{equation}
\label{bound:ddot Y 2}
 \Vert \ddot Y^\gd_t\Vert \, \leq\,  C(\gd+\gd \Vert \dot Y^\gd_t\Vert)\, .
\end{equation}
Suppose that $\sup_{t\leq \tau^\gd_{2}} \Vert \dot Y^\gd_t\Vert$ is reached at $\Vert\dot Y^\gd_{t_0}\Vert$. For $\gd$ small enough, for $s\in[t_0-1,t_0]$ we have
 \begin{equation}
  \label{ineq:bound dot Y below 2}
 \Vert\dot Y^\gd_s\Vert\, \geq\, \frac12\Vert \dot Y^\gd_{t_0}\Vert-C\gd\, .
\end{equation}
Proceeding as in Lemma \ref{lem:dot Y 1/2}, we can bound the cost of the path $Y^\gd$ on the time interval $[t_0-1,t_0]$ from below by
\begin{equation}
\label{ineq:cost Y}
 \frac14 \int_{t_0-1}^{t_0}\Vert \dot Y^\gd_t\Vert^2\dd t-\frac12\int_{t_0-1}^{t_0}\Vert \nabla V[Y^\gd_t] -\gd G[Y^\gd_t]\Vert^2\dd t\, \geq \, \frac14 \left(\frac12\Vert \dot Y^\gd_{t_0}\Vert-C\gd\right)^2-C\gd^2\, .
\end{equation}
On the other hand, we consider the path defined on the time interval $[0,\Vert  \dot Y^\gd_t\Vert/\gd]$ by
\begin{equation}
 Z_t\, =\, Y^\gd_{t_0-1+\gd t/\Vert \dot Y^\gd_t\Vert}\, ,
\end{equation}
which connects the points $Y^\gd_{t_0-1}$ and $Y^\gd_{t_0}$. Using Lemma \ref{lem:dist gd M} we can bound the cost of this path:
\begin{equation}
\label{ineq:cost Z}
 \int_{0}^{\Vert  \dot Y^\gd_t\Vert/\gd}\Vert \dot Z_t+\nabla V[Z_t]-\gd G[Z_t]\Vert^2\dd t\, \leq\, C\frac{\Vert  \dot Y^\gd_t\Vert}{\gd}\left(\sup_{t\in [0,\Vert  \dot Y^\gd_t\Vert/\gd]}\Vert \dot Z_t\Vert^2+\gd^2\right)\, \leq\, C\gd \Vert  \dot Y^\gd_t\Vert\, .
\end{equation}
So if $\Vert  \dot Y^\gd_t\Vert$ is too large, that is $\Vert  \dot Y^\gd_t\Vert\geq C_4\gd$ for $C_4$ large enough, we get a contradiction between \eqref{ineq:cost Y} and \eqref{ineq:cost Z}: by distending time and replacing $Y^\gd_{|(t_0-1,t_0)}$ by $Z$ we can create better path than $Y^\gd$ for $B^\gd$.
\end{proof}

\medskip

\begin{lemma}
\label{lem:dist gd2 M}
There exists $C_5>0$ such that for all $\gd$ small enough, if we define
\begin{equation}
 \tau^\gd_{1}\, :=\, \inf \{t<0, \Vert h^\gd_t-\phi_\gd(\theta^\gd_t)\Vert\geq C_5\gd^2\}\, ,
\end{equation}
then $\Vert h^\gd_t-\phi_\gd(\theta^\gd_t)\Vert \geq C_5\gd^2$ for all $\tau^\gd_{1}\leq t\leq 0$ and if $\tau^\gd_{1}<0$, $\Vert h^\gd_t-\phi_\gd(\theta^\gd_t)\Vert$ increases strictly on $(\tau^\gd_{1},0]$.
Moreover there exists $C_6>0$ such that
\begin{equation}
 \tau^\gd_{1} \, \geq\, -C_6  |\log\gd|\, .
\end{equation}
\end{lemma}

\medskip

\begin{proof}
We proceed as for Lemma \ref{lem:dist gd M}. Using Lemma \ref{lem:first order phi}, Lemma \ref{lem:dist gd M}, Lemma \ref{lem:dot Y gd} and \eqref{eq:zero grad} we get for $t\leq \tau^\gd_{2}$ the following first order expansion of \eqref{eq:EL}:
\begin{equation}
\label{eq:EL first order 2}
 \ddot Y^\gd_t\, =\, \Big(H[q(\theta^\gd_t)]\Big)^2(h^\gd_t-\phi_\gd(\theta^\gd_t))+O(\gd^2)\, .
\end{equation}
This time we have (recall \eqref{eq:def alpha})
\begin{equation}
\label{eq:ddot ga first order}
 \ddot \ga^\gd_t\, =\,2\Vert \dot h^\gd_t -\dot\theta^\gd_t\phi_\gd'(\theta^\gd_t)\Vert^2+ 2\Vert H[q(\theta^\gd_t)](h^\gd_t-\phi_\gd(\theta^\gd_t))\Vert^2+O(\gd^2\Vert h^\gd_t-\phi_\gd(\theta^\gd_t)\Vert)\, .
\end{equation}
So (recall \eqref{eq:spect gap hess}) there exists $C>0$ such that
\begin{equation}
 \ddot \ga^\gd_t\, \geq\, 2\gl^2 \ga^\gd_t-C\gd^2 (\ga^\gd_t)^{1/2}\, .
\end{equation}
In particular, for $ (\ga^\gd_t)^{1/2}\geq C_5\gd^2$ where $C_5:=C/\gl^2$, we have
\begin{equation}
 \ddot \ga^\gd_t\, \geq\,\gl^2 \ga^\gd_t\, .
\end{equation}
We deduce that if $\Vert h^\gd_t-\phi_\gd(\theta^\gd_t)\Vert\geq C_5\gd^2$ and $\dot\ga^\gd_t\geq 0$ for a $t<0$, then $\ga^\gd_s$ is strictly increasing for $s>t$. Moreover, it takes at most a time of order $|\log\gd|$ for $\ga^\gd$ to reach $C_0\gd^{1/2}$.

\end{proof}

\medskip

\begin{lemma}
\label{lem:dot h 2}
There exists $C_7>0$ such that for all $\gd$ small enough 
\begin{equation}
 \sup_{t\leq \tau^\gd_{1}}\Vert \dot h^\gd_t\Vert\, \leq\, C_7\gd^2\, .
\end{equation}
\end{lemma}

\medskip

\begin{proof}
Suppose that $\sup_{t\leq \tau^\gd_{1}} \Vert \dot h^\gd_t\Vert$ is reached at $\Vert\dot h^\gd_{t_0}\Vert$. Then the mean value Theorem implies
\begin{multline}
 \Big\Vert h^\gd_{t_0-1}-\phi_\gd(\theta^\gd_{t_0-1})-\Big( h^\gd_{t_0}- \phi_\gd(\theta^\gd_t)+ \dot h^\gd_{t_0}-\dot\theta^\gd_{t_0}\phi_\gd'(\theta^\gd_{t_0})\Big)\Big\Vert \\ \leq\, \sup_{t\in[t_0-1,t_0]}\Vert \ddot h^\gd_t-\ddot\theta^\gd_t\phi_\gd'(\theta^\gd_t)-(\dot\theta^\gd_t)^2\phi_\gd''(\theta^\gd_t)\Vert\, .
\end{multline}
Theorem \ref{lem:phi Mgd}, Lemma \ref{lem:dot Y gd}, Lemma \ref{lem:dist gd2 M} and \eqref{eq:EL first order 2} imply that the right hand side is of order $\gd^2$. But the same arguments imply that all the right hand side with $\dot h^\gd_{t_0}$ taken away is also of order $\gd^2$. We deduce that $\dot h^\gd_{t_0}$ must also be of order $\gd^2$.
\end{proof}



\medskip

\begin{lemma}
\label{lem:exit}
For $\gd$ small enough, the projection $p_\gd(B^\gd)$ on $M^\gd$ coincides either with $\varphi_{A^\gd}-\Delta_1$ or with $\varphi_{A^\gd}+\Delta_2$.
Moreover, there exist $C_8$, $C_9$ and $C_{10}$ such that for $\gd$ small enough and for all $\tau^\gd_1\leq t\leq 0$
\begin{equation}
\label{eq:bound dot theta exit}
 \vert \dot\theta^\gd_t\vert \, \leq\, C_8 \gd|\log\gd|\, ,
\end{equation}
\begin{equation}
\label{eq:bound h exit}
 \Vert h^\gd_t \Vert\, \leq\, C_9 \gd |\log\gd|\, ,
\end{equation}
and
\begin{equation}
\label{eq:bound dot Y exit}
 \Vert \dot Y^\gd_t\Vert\, \leq\, C_{10} \gd|\log\gd|\, .
\end{equation}

\end{lemma}

\medskip
\begin{rem}
\label{rem:exit}
 This Lemma implies in particular that $\sup_{t\leq 0}\dist(Y^\gd,M^\gd)=O(\gd|\log\gd|)$, and thus recalling Lemma \ref{lem:prelim exit Y} that $Y^\gd$ does not exit prematurely  $U^\gd$.
\end{rem}

\medskip

\begin{proof}
Isolating the $\gd G$ term in $I_{\gd,\tau^\gd_1}(Y^\gd)$, we get
\begin{equation}
 I_{\gd,\tau^\gd_1}(Y^\gd)\, \geq\, \frac14\int_{\tau^\gd_1}^0 \Vert \dot Y^\gd_t +\nabla V[Y^\gd_t]\Vert^2\dd t-\frac{\gd^2}{2} \int_{\tau^\gd_1}^0 \Vert G[Y^\gd_t]\Vert^2\dd t\, .
\end{equation}
Thanks to Lemma \ref{lem:dist gd2 M}, we know that the last term of the right hand side is of order $\gd^2|\log\gd|$, so inserting the reversed time dynamics (using the identity $\Vert u-v\Vert^2=\Vert u+v\Vert^2-4\langle u,v\rangle$), we obtain for a $C>0$
\begin{equation}
\label{ineq:I gd tau1}
 I_{\gd,\tau^\gd_1}(Y^\gd)\, \geq\,  \frac14\int_{\tau^\gd_2}^0 \Vert \dot Y^\gd_t -\nabla V[Y^\gd_t]\Vert^2\dd t+ V[B^\gd]-V[Y^\gd_{\tau^\gd_1}]-C\gd^2|\log\gd|\, .
\end{equation}
Since $Y^\gd_{\tau^\gd_1}$ is located at a distance of order $\gd$ from $M$, $V[Y^\gd_{\tau^\gd_1}]$ is of order $\gd^2$. Thus \eqref{ineq:I gd tau1} implies 
\begin{equation}
\label{ineq:I gd tau1 2}
 I_{\gd,\tau^\gd_1}(Y^\gd)\, \geq\, V[B^\gd]-C\gd^2|\log\gd|\, .
\end{equation}
Taking \eqref{ineq:first bound W} and \eqref{eq:choice C0} into account, this implies the first assertion of the Lemma for $\gd$ small enough.

Lemma \ref{lem:p} implies for all $t\leq 0$
\begin{equation}
 \dot \theta^\gd_t\, =\, \frac{1}{1-\langle h^\gd_t,q''(\theta^\gd_t)\rangle}\langle \dot Y^\gd_t,q'(\theta^\gd_t)\rangle\, .
\end{equation}
So a derivation in time gives
\begin{multline}
\label{eq:ddot theta}
 \ddot \theta^\gd_t\, =\,  \frac{1}{1-\langle h^\gd_t,q''(\theta^\gd_t)\rangle}\langle \ddot Y^\gd_t,q'(\theta^\gd_t)\rangle+ \frac{\dot \theta^\gd_t}{1-\langle h^\gd_t,q''(\theta^\gd_t)\rangle}\langle \dot Y^\gd_t,q''(\theta^\gd_t)\rangle\\ + \frac{\langle \dot h^\gd_t,q''(\theta^\gd_t)\rangle+\dot\theta^\gd_t\langle h^\gd_t,q'''(\theta^\gd_t)\rangle}{(1-\langle h^\gd_t,q''(\theta^\gd_t)\rangle)^2}\langle \dot Y^\gd_t,q'(\theta^\gd_t)\rangle\, .
\end{multline}
Lemma \ref{lem:dot Y 1/2} implies that the two last terms of the right hand side are of order $\gd$. Moreover \eqref{eq:EL first order} and \eqref{eq:zero hess} imply that the first term of the right hand side is also of order $\gd$. We deduce that theres exists $C>0$ such that for all $\tau^\gd_1\leq t\leq 0$
\begin{equation}
 |\dot \theta^\gd_t-\dot\theta^\gd_{\tau^\gd_1}|\, \leq\, C\gd|\log\gd|\, .
\end{equation}
To deduce \eqref{eq:bound dot theta exit}, we just have to remind Lemma \ref{lem:dot Y gd}, and more precisely that $\dot\theta^\gd_{\tau^\gd_1}$ is of order $\gd$.

\medskip

To get the two last results of Lemma \ref{lem:exit}, we find an upper bound for $I_{\gd,\tau^\gd_1}(Y^\gd)$ by studying the cost of a specifically chosen path. We consider a path $Z_t$ defined on the time interval $[T,0]$, starting from $Y^\gd_{\tau^\gd_1}$, linking in a linear way $Y^\gd_{\tau^\gd_1}$ to $\tilde q_\gd(\theta^\gd_{\tau^\gd_1})$, and then following the curve $M^\gd$ to exit $U^\gd$ at the point $q_\gd(\varphi_{B^\gd})$ (we know that it is possible to exit $U^\gd$ this way, thanks to the first assertion of the Lemma). Since $Y^\gd$ is an optimal path for the exit from $U^\gd$, it is clear that
$ \int_T^{0}\Vert \dot Z_t+\nabla V[Z_t]-\gd G[Z_t]\Vert^2\dd t $ is an upper bound for $I_{\gd,\tau^\gd_1}(Y^\gd)$. More precisely we define for $T\leq t\leq T+\gd$:
\begin{equation}
 Z_t\, :=\, Y^\gd_{\tau^\gd_1}+ \frac{t-T}{\gd}\Big(\tilde q_\gd(\theta^\gd_{\tau^\gd_1})-Y^\gd_{\tau^\gd_1}\Big)\, =\, Y^\gd_{\tau^\gd_1}-\frac{t-T}{\gd}h^\gd_{\tau^\gd_1}\, .
\end{equation}
In this case, using \eqref{eq:zero grad} and the fact that $h^\gd_{\tau^\gd_1}$ is of order $\gd^2$, we get the bound
\begin{multline}
  \int_T^{T+\gd}\Vert \dot Z_t+\nabla V[Z_t]-\gd G[Z_t]\Vert^2\dd t\, \leq\, 2\int_T^{T+\gd}\Vert \dot Z_t\Vert^2\dd t+2\int_T^{T+\gd}\Vert \nabla V[Z_t]-\gd G[Z_t]\Vert^2\dd t
  \\ =\, O(\gd^3)
\end{multline}
On the order hand, if we define $\varphi_{\tau^\gd_1}$ the phase satisfying $q_\gd(\varphi_{\tau^\gd_1})=\tilde q_\gd(\theta^\gd_{\tau^\gd_1})$, since for $t\geq T+\gd$ $Z_t$ follows the curve $M^\gd$, we can make $\int_{T+\gd}^{0}\Vert \dot Z_t+\nabla V[Z_t]-\gd G[Z_t]\Vert^2\dd t$ as close to $\int_{\varphi_{\tau^\gd_1}}^{\varphi_{B^\gd}} b_\gd(\varphi)\dd \varphi$ as we want. But \eqref{ineq:I gd tau1} implies that $\theta^\gd_0-\theta^\gd_{\tau^\gd_1}$ is of order $\gd |\log\gd|^2$, and Lemma \ref{lem:p close to tilde p} implies that $\varphi_{B^\gd}-\varphi_{\tau^\gd_1}$ is of the same order. We deduce that $\int_{\varphi_{\tau^\gd_1}}^{\varphi_{B^\gd}} b_\gd(\varphi)\dd \varphi$ is of order $\gd^2 |\log(\gd)|^2$, and thus that $I_{\gd,\tau^\gd_1}(Y^\gd)$ is also at most of order 
$\gd^2 |\log(\gd)|^2$.

Recalling \eqref{ineq:I gd tau1 2}, \eqref{eq:zero grad} and \eqref{eq:zero hess}, we deduce that $\dist(B^\gd,M^\gd)$ is at most of order $\gd|\log\gd|$, and thus we get \eqref{eq:bound h exit}, since $M^\gd$ is located at a distance of order $\gd$ from $M$. The proof of \eqref{eq:bound dot Y exit} is similar to the proof of Lemma \ref{lem:dot Y gd}, keeping \eqref{eq:bound h exit} in mind.
\end{proof}

\end{subsection}

\medskip

\begin{subsection}{Proof of Theorem \lowercase{\ref{th:exp w}}}

We will show that $I_{\gd,-\infty}(Y^\gd)$ can be well approximated by $I_{\gd,\tau^\gd_0}(Y^\gd)$ with a good choice of $\tau^\gd_0\leq \tau^\gd_1$, and then we will study the cost of the path $Y^\gd$ on the time intervals $[\tau^\gd_0,\tau^\gd_1]$ (when $\tau^\gd_0<\tau^\gd_1$) and $[\tau^\gd_1,0]$.

\medskip

Recalling \eqref{eq:zero grad}, Lemma \ref{lem:first order phi}, Lemma \ref{lem:dist gd2 M} and Lemma \ref{lem:dot Y gd}, we can expand \eqref{eq:EL} in the following way for all $t\leq \tau^\gd_1$: 
\begin{multline}
 \ddot Y^\gd_t\, =\, \Bigg( H[q(\theta^\gd_t)] + D^3V[q(\theta^\gd_t)](h^\gd_t,\ldotp,\ldotp)-\gd  DG^\dagger[q(\theta^\gd_t)] +O(\gd^2)\Bigg)
\\ \times\Bigg(H[q(\theta^\gd_t)]h^\gd_t-\gd \Big(G[q(\theta^\gd_t)]-\langle G[q(\theta^\gd_t)],q'(\theta^\gd_t)\rangle q'(\theta^\gd_t)\Big)+D^3V[q(\theta^\gd_t)](h^\gd_t,h^\gd_t,\ldotp)
\\ -\gd \langle G[q(\theta^\gd_t)],q'(\theta^\gd_t)\rangle q'(\theta^\gd_t)-\gd DG[q(\theta^\gd_t)]h^\gd_t+O(\gd^3)\Bigg)
 \\ +\gd \Big( DG[q(\theta^\gd_t)]-  DG^\dagger[q(\theta^\gd_t)]\Big) \dot Y^\gd_t +O(\gd^3)\, ,
\end{multline}
and (recall \eqref{eq:zero hess}) we expand the product in the following way:
\begin{multline}
 \ddot Y^\gd_t\, =\, H[q(\theta^\gd_t)]\Bigg(H[q(\theta^\gd_t)]h^\gd_t-\gd \Big(G[q(\theta^\gd_t)]-\langle G[q(\theta^\gd_t)],q'(\theta^\gd_t)\rangle q'(\theta^\gd_t)\Big)\Bigg)
\\  +H[q(\theta^\gd_t)]D^3V[q(\theta^\gd_t)](h^\gd_t,h^\gd_t,\ldotp)+\gd H[q(\theta^\gd_t)]DG[q(\theta^\gd_t)]h^\gd_t
 \\  +\gd \langle G[q(\theta^\gd_t)],q'(\theta^\gd_t)\rangle D^3V[q(\theta^\gd_t)](h^\gd_t,q'(\theta^\gd_t),\ldotp)+\gd^2 \langle G[q(\theta^\gd_t)],q'(\theta^\gd_t)\rangle DG^\dagger[q(\theta^\gd_t)] q'(\theta^\gd_t)
 \\ +\gd \Big( DG[q(\theta^\gd_t)]- DG^\dagger[q(\theta^\gd_t)]\Big) \dot Y^\gd_t +O(\gd^3)
\end{multline}
Using again \eqref{eq:zero hess} we get the projection
\begin{multline}
\label{eq:proj ddot Y}
\langle \ddot Y^\gd_t,q'(\theta^\gd_t)\rangle\, =\, \gd \Big\langle  G[q(\theta^\gd_t)],q'(\theta^\gd_t)\Big\rangle D^3V[q(\theta^\gd_t)](h^\gd_t,q'(\theta^\gd_t),q'(\theta^\gd_t))       \\+ \gd^2 \Big\langle  G[q(\theta^\gd_t)],q'(\theta^\gd_t)\Big\rangle \Big\langle  DG^\dagger[q(\theta^\gd_t)] q'(\theta^\gd_t),q'(\theta^\gd_t)\Big\rangle
\\+\gd \Big\langle\Big( DG[q(\theta^\gd_t)]-  DG^\dagger[q(\theta^\gd_t)]\Big) \dot Y^\gd_t,q'(\theta^\gd_t)\Big\rangle +O(\gd^3)\, .
\end{multline}
But Lemma \ref{lem:dot h 2} implies the following first order expansion of $\dot Y^\gd_t$ for $t\leq\tau^\gd_1$:
\begin{equation}
 \dot Y^\gd_t\, =\, \dot \theta^\gd_t q'(\theta^\gd_t)+O(\gd^2)\, ,
\end{equation}
and thus the first term of the last line in \eqref{eq:proj ddot Y} is in fact 
\begin{multline}
 \Big\langle\Big( DG[q(\theta^\gd_t)]-  DG^\dagger[q(\theta^\gd_t)]\Big) \dot Y^\gd_t,q'(\theta^\gd_t)\Big\rangle\, =\, \dot\theta^\gd_t\langle DG[q(\theta^\gd_t)]q'(\theta^\gd_t),q'(\theta^\gd_t)\rangle \\-\dot\theta^\gd_t\langle q'(\theta^\gd_t),DG[q(\theta^\gd_t)]q'(\theta^\gd_t)\rangle +O(\gd^3)\, =\, O(\gd^3).
\end{multline}
Applying \eqref{eq:link V3 H} to \eqref{eq:proj ddot Y}, we get
\begin{multline}
 \langle \ddot Y^\gd_t,q'(\theta^\gd_t)\rangle\, =\, \gd^2 \Big(\langle DG[q(\theta^\gd_t) q'(\theta^\gd_t),q'(\theta^\gd_t)\rangle +\langle G[q(\theta^\gd_t)],q''(\theta^\gd_t)\rangle \Big)\langle G[q(\theta^\gd_t)],q'(\theta^\gd_t)\rangle+O(\gd^3)
 \\ =\, b'(\theta^\gd_t)b(\theta^\gd_t)+O(\gd^3)\, .
\end{multline}
Now recalling \eqref{eq:ddot theta}, Lemma \ref{lem:dist gd2 M}, Lemma \ref{lem:dot Y gd} and Lemma \ref{lem:dot h 2} we get the same expansion for $\ddot \theta^\gd_t$: for all $t\leq \tau^\gd_1$ we have
\begin{equation}
\label{eq:ddot theta b}
\ddot \theta^\gd_t\, =\,  b'(\theta^\gd_t)b(\theta^\gd_t)+O(\gd^3) \, .
\end{equation}
Remind that we have supposed that the dynamics $\dot\theta = b(\theta)$ admits an hyperbolic stable fixed point $\theta^0$, and remind  (recall Lemma \ref{lem:p close to tilde p} and Lemma \ref{lem:tilde b}) that $\theta_{A^\gd}:=p(A^\gd)$ is located at a distance of order $\gd$ from $\theta^0$. So when $\theta^\gd_t$ is close to $\theta_{A^\gd}$, though at a distance greater than $\gd$, the leading term in \eqref{eq:ddot theta b} is up to a constant factor equivalent to $\gd^2(\theta^\gd_t-\theta_{A^\gd})$. If we suppose (without loss of generality) that $p_\gd(B^\gd)$ is located on the side of the increasing $\theta$ with respect to $\theta_{A^\gd}$, and we define $\tau^\gd_0$ the first time $t\leq\tau^\gd_1$ such that $\theta^\gd_t-\theta^\gd_{A^\gd}=\gd|\log\gd|$, then $\theta^\gd_t$ is exponentially increasing (at rate $\exp(C\gd t)$) for $t\geq \tau^\gd_0$, at least until it reaches $1/|\log\gd|$ at a time $t_1$. So $t_1-\tau^\gd_0=O(|\log\gd|/\gd)
$. 
We now want to bound the difference $\tau^\gd_1-\tau^\gd_0$. If $t_1<\tau^\gd_1$, then for $t_1\leq t\leq \tau^\gd_1$, by multiplying each side
of \eqref{eq:ddot theta b} by $\dot\theta^\gd_t$ and integrating in time, we get (recall that $\dot\theta^\gd_t=O(\gd)$ for $t\leq \tau^\gd_1$)
\begin{equation}
\label{eq:dot theta squared}
 (\dot\theta^\gd_t)^2-(\dot\theta^\gd_{\tau^\gd_0})^2\, =\, b^2(\theta^\gd_t)-b^2(\theta^\gd_{\tau^\gd_0})+O(\gd^4(t-\tau^\gd_0))\, .
\end{equation}
By construction $b(\theta^\gd_{\tau^\gd_0})=O(\gd^2|\log\gd|)$ and $b(\theta_{t_1})$ is up to a constant factor equivalent to $\gd/|\log\gd|$. Moreover, if we denote $\theta_{B^\gd}:=p(B^\gd)$, since $[\theta^\gd_{t_1},\theta_{B^\gd}]$  is strictly included in the domain of attraction of $\theta^0$ for $b$, for $\gd$ small enough $|b(\theta)|$ is greater than $\gd/|\log\gd|$ for $\theta\in [\theta^\gd_{t_1},\theta_{B^\gd}]$. Thus $\theta^\gd_t$ keeps increasing for $t\geq t_1$ and  \eqref{eq:dot theta squared} implies
\begin{equation}
 \dot\theta^\gd_t\, \geq\, |b(\theta^\gd_t)|+O(\gd^3|\log\gd|,\gd^3|\log\gd| (t-\tau^\gd_0))\, .
\end{equation}
Dividing by $|b(\theta^\gd_t)|$ and integrating in time between $t_1$ and $\tau^\gd_1$, we see that $\tau^\gd_1-t_1=O(|\log\gd|/\gd)$. We deduce
\begin{equation}
\label{ineq:bound tau1 minus tau0}
 \tau^\gd_1-\tau^\gd_0\, =\, O(|\log\gd|/\gd)\, .
\end{equation}

\medskip

We easily see that the cost of the path $Y^\gd$ restricted to the times smaller than $\tau^\gd_0$ is negligible. Since $Y^\gd$ is an optimal path, this cost is exactly $W(A^\gd,Y^\gd_{\tau^\gd_0})$, and this quasipotential is of order $\gd^3(\log\gd)^2$. Indeed, reaching $\tilde q_\gd(\theta^\gd_{\tau^\gd_0})$ following $M^\gd$ costs $\int_{\varphi_{A^\gd}}^{\varphi_{\tau^\gd_0}}b_\gd(\varphi)\dd\varphi$ where $\varphi_{\tau^\gd_0}$ satisfies $q(\varphi_{\tau^\gd_0})=\tilde q_\gd(\theta^\gd_{\tau^\gd_0})$, and (recall Lemma \ref{lem:p close to tilde p} and the definition of $\tau^\gd_0$) $\varphi_{\tau^\gd_0}-\varphi_{A^\gd}\sim \gd|\log\gd|$. Since $b_\gd(\varphi_{A^\gd})=0$ and $b_\gd(\phi)=O(\gd)$, we indeed get $\int_{\varphi_{A^\gd}}^{\varphi_{\tau^\gd_0}}b_\gd(\varphi)\dd\varphi=O(\gd^3|\log\gd|^2)$. Moreover it is possible to link linearly $q_\gd(\theta^\gd_{\tau^\gd_0})$ to $Y^\gd_{\tau^\gd_0}$ at cost $\gd^3$ (see the proof of Lemma \ref{lem:exit}).

\medskip

We now study the cost of the $Y^\gd_t$ on the time interval $[\tau^\gd_0,\tau^\gd_1]$. It is possible to get a lower bound by projecting on the tangent space of $M$:
\begin{equation}
\label{eq:I projected}
 \frac12\int_{\tau^\gd_0}^{\tau^\gd_1} \Vert \dot Y^\gd_t+\nabla V[Y^\gd_t]-\gd G[Y^\gd_t]\Vert^2\dd t\, \geq\, \frac12\int_{\tau^\gd_0}^{\tau^\gd_1} |\langle \dot Y^\gd_t+\nabla V[Y^\gd_t]-\gd G[Y^\gd_t],q'(\theta^\gd_t)\rangle|^2\dd t\, .
\end{equation}
It is now enough to simply expand the integrand. Lemma \ref{lem:p}, Lemma \ref{lem:dist gd2 M} and Lemma \ref{lem:dot Y gd} imply
\begin{equation}
\label{eq:proj dot Y}
 \langle \dot Y^\gd_t,q'(\theta^\gd_t)\rangle \, =\, \dot\theta^\gd_t(1 -\langle h^\gd_t,q''(\theta^\gd_t)\rangle)\, =\, \dot\theta^\gd_t(1 -\langle \phi_\gd(\theta^\gd_t),q''(\theta^\gd_t)\rangle)+O(\gd^3)
\end{equation}
and (using moreover \eqref{eq:zero grad} and \eqref{eq:zero hess})
\begin{multline}
\label{eq:proj nabla V G}
 \langle\nabla V[Y^\gd_t]-\gd G[Y^\gd_t], q'(\theta^\gd_t)\rangle\, =\,  D^3V[q(\theta^\gd_t)](h^\gd_t,h^\gd_t,q'(\theta^\gd_t)) \\-\gd \langle G[q(\theta^\gd_t)]+DG[q(\theta^\gd_t)]h^\gd_t,q'(\theta^\gd_t)\rangle  +O(\gd^2) \\
 =\, D^3V[q(\theta^\gd_t](\phi_\gd(\theta^\gd_t),\phi_\gd(\theta^\gd_t),q'(\theta^\gd_t)) \\-\gd\langle G[q(\theta^\gd_t)]+DG[q(\theta^\gd_t)]\phi_\gd(\theta^\gd_t),q'(\theta^\gd_t)\rangle  +O(\gd^2)\, .
\end{multline}
Remark that Theorem \ref{lem:phi Mgd}, \eqref{eq:zero grad} and \eqref{eq:zero hess} also imply
\begin{multline}
 \langle \nabla V[\tilde q_\gd(\theta^\gd_t)]-\gd G[\tilde q_\gd(\theta^\gd_t)], q'(\theta^\gd_t)\rangle\, =\, D^3V[q(\theta^\gd_t)](\phi_\gd(\theta^\gd_t),\phi_\gd(\theta^\gd_t),q'(\theta^\gd_t)) \\-\gd \langle G[q(\theta^\gd_t)]+DG[q(\theta^\gd_t)]\phi_\gd(\theta^\gd_t),q'(\theta^\gd_t)\rangle  +O(\gd^2)\, ,
\end{multline}
so
\begin{equation}
  \langle\nabla V[Y^\gd_t]-\gd G[Y^\gd_t], q'(\theta^\gd_t)\rangle\, =\, \langle \nabla V[\tilde q_\gd(\theta^\gd_t)]-\gd G[\tilde q_\gd(\theta^\gd_t)], q'(\theta^\gd_t)\rangle+O(\gd^2)\, .
\end{equation}
Since $M^\gd$ is invariant under \eqref{eq:perturbed1 determ}, $\nabla V[\tilde q_\gd(\theta^\gd_t)]-\gd G[\tilde q_\gd(\theta^\gd_t)]$ and $\tilde q'_\gd(\theta^\gd_t)$ are collinear vectors, and thus 
\begin{equation}
\label{eq:nabla V G proj qgd}
 \langle \nabla V[Y^\gd_t]-\gd G[Y^\gd_t], q'(\theta^\gd_t)\rangle \, =\,\langle\nabla V[\tilde q_\gd(\theta^\gd_t)]-\gd G[\tilde q_\gd(\theta^\gd_t)], \tilde q'_\gd(\theta^\gd_t)\rangle \frac{\langle \tilde q'_\gd(\theta^\gd_t),q'(\theta^\gd_t)\rangle}{\Vert \tilde q'_\gd(\theta^\gd_t)\Vert^2 }+O(\gd^2)\, .
\end{equation}
Now $\langle \phi_\gd(\theta),q'(\theta)\rangle=0$ implies $\langle \phi_\gd'(\theta), q'(\theta)\rangle =-\langle \phi_\gd(\theta), q''(\theta)\rangle$ and thus using Theorem \ref{lem:phi Mgd} we get
\begin{equation}
\label{eq:norm q phi}
 \Vert \tilde q'_\gd(\theta^\gd_t)\Vert\, =\, (1+2\langle \phi'_\gd(\theta^\gd_t),q'(\theta^\gd_t)\rangle+O(\gd^2))^{1/2}\,  =\, 1-\langle \phi_\gd(\theta^\gd_t),q''(\theta^\gd_t)\rangle+O(\gd^2)\, ,
\end{equation}
and similarly
\begin{equation}
\label{eq:q phi scalair q}
 \langle \tilde q'_\gd(\theta^\gd_t),q'(\theta^\gd_t)\rangle \, =\, 1-\langle \phi_\gd(\theta^\gd_t),q''(\theta^\gd_t)\rangle+O(\gd^2)
 \, =\, \Vert \tilde q'_\gd(\theta^\gd_t)\Vert+O(\gd^2)\, .
\end{equation}
So in view of \eqref{eq:proj dot Y}, \eqref{eq:nabla V G proj qgd}, \eqref{eq:norm q phi}, \eqref{eq:q phi scalair q} and \eqref{ineq:bound tau1 minus tau0}, \eqref{eq:I projected} becomes
\begin{multline}
 \frac12\int_{\tau^\gd_0}^{\tau^\gd_1} \Vert \dot Y^\gd_t+\nabla V[Y^\gd_t]-\gd G[Y^\gd_t]\Vert^2\dd t \\ \geq\, \frac12 \int_{\tau^\gd_0}^{\tau^\gd_1} \bigg|\dot  \theta^\gd_t  \Vert \tilde q'_\gd(\theta^\gd_t)\Vert -  \bigg\langle \nabla V[\tilde q_\gd(\theta^\gd_t)]-\gd G[\tilde q_\gd(\theta^\gd_t)], \frac{\tilde q'_\gd(\theta^\gd_t)}{\Vert  \tilde q'_\gd(\theta^\gd_t)\Vert}\bigg\rangle\bigg|^2  
  \dd t   +O(\gd^{3}|\log\gd|)\, .
\end{multline}

\medskip

We proceed similarly on the time interval $[\tau^\gd_1,0]$: we also use the lower bound
\begin{equation}
\label{eq:I projected 2}
 \frac12\int_{\tau^\gd_1}^{0} \Vert \dot Y^\gd_t+\nabla V[Y^\gd_t]-\gd G[Y^\gd_t]\Vert^2\dd t\, \geq\, \frac12\int_{\tau^\gd_1}^{0} |\langle \dot Y^\gd_t+\nabla V[Y^\gd_t]-\gd G[Y^\gd_t],q'(\theta^\gd_t)\rangle|^2\dd t\, .
\end{equation}
This time, using Lemma \ref{lem:exit}, \eqref{eq:norm q phi} and \eqref{eq:q phi scalair q} we get
\begin{equation}
 \dot Y^\gd_t\, =\, \dot\theta^\gd_t \Vert \tilde q'_\gd(\theta^\gd_t)\Vert+O(\gd^2(\log\gd)^2)\, 
\end{equation}
and
\begin{equation}
 \langle\nabla V[Y^\gd_t]-\gd G[Y^\gd_t], q'(\theta^\gd_t)\rangle\, =\, \bigg\langle \nabla V[\tilde q_\gd(\theta^\gd_t)]-\gd G[\tilde q_\gd(\theta^\gd_t)], \frac{\tilde q'_\gd(\theta^\gd_t)}{\Vert  \tilde q'_\gd(\theta^\gd_t)\Vert}\bigg\rangle+O(\gd^2(\log\gd)^2)\, .
\end{equation}
Since $\tau^\gd_1=O(|\log\gd|)$, we deduce
\begin{multline}
 \frac12\int_{\tau^\gd_1}^{0} \Vert \dot Y^\gd_t+\nabla V[Y^\gd_t]-\gd G[Y^\gd_t]\Vert^2\dd t \\ \geq\, \frac12 \int_{\tau^\gd_1}^{0} \bigg|\dot  \theta^\gd_t  \Vert \tilde q'_\gd(\theta^\gd_t)\Vert -  \bigg\langle \nabla V[\tilde q_\gd(\theta^\gd_t)]-\gd G[\tilde q_\gd(\theta^\gd_t)], \frac{\tilde q'_\gd(\theta^\gd_t)}{\Vert  \tilde q'_\gd(\theta^\gd_t)\Vert}\bigg\rangle\bigg|^2  
  \dd t   +O(\gd^{3}|\log\gd|^3)\, .
\end{multline}

In conclusion we have proved
\begin{equation}
\label{eq:I approx}
I_{\gd,-\infty}^{A^\gd}(Y^\gd)\, \geq\, \frac12 \int_{\tau^\gd_0}^{0} \bigg|\dot  \theta^\gd_t  \Vert \tilde q'_\gd(\theta^\gd_t)\Vert -  \bigg\langle \nabla V[\tilde q_\gd(\theta^\gd_t)]-\gd G[\tilde q_\gd(\theta^\gd_t)], \frac{\tilde q'_\gd(\theta^\gd_t)}{\Vert  \tilde q'_\gd(\theta^\gd_t)\Vert}\bigg\rangle\bigg|^2  
  \dd t   +O(\gd^{3}|\log\gd|^3)\, .
\end{equation}
If we define $\varphi^\gd_t:=\tilde p_\gd(\tilde q_\gd(\theta^\gd_t))$ for all $\tau^\gd_0\leq t\leq 0$ (i.e $q_\gd(\varphi^\gd_t)=\tilde q_\gd(\theta^\gd_t)$), then \eqref{eq:I approx} becomes
\begin{equation}
\label{eq:last approx I}
 I_{\gd,-\infty}^{A^\gd}(Y^\gd)\, \geq\, \frac12 \int_{\tau^\gd_0}^{0} |\dot  \varphi^\gd_t  - b_\gd(\varphi^\gd_t)|^2  
  \dd t   +O(\gd^{3}|\log\gd|^3)\, .
\end{equation}
We stress that $\varphi^\gd_t$ is not exactly $p_\gd(Y^\gd_t)$, but Lemma \ref{lem:p close to tilde p} and Lemma \ref{lem:exit} ensure that the induced error is negligible. We have indeed $\varphi^\gd_0-\varphi_{B^\gd}=O(\gd^2|\log\gd|)$. On the other hand $\varphi^\gd$ does not begin at the point $\varphi_{A^\gd}$, but remind that $\varphi^\gd_{\tau^\gd_0}-\varphi_{A^\gd}=O(\gd|\log\gd|)$, so
$W^{red}(\varphi_{A^\gd},\varphi^\gd_{\tau^\gd_0})=O(\gd^3|\log\gd|^2)$. These two observations imply
\begin{multline}
 \inf\left\{\frac12 \int_{\tau^\gd_0}^{0} |\dot  \varphi_t  - b_\gd(\varphi_t)|^2  \dd t:\, \varphi \text{ is }C^2,\,  \varphi_{\tau^\gd_0}=\varphi^\gd_{\tau^\gd_0}\text{ and }\varphi_0=  \varphi^\gd_0\right\} \\ \geq\, W^{red}_\gd( \varphi_{A^\gd},\varphi_{B^\gd})+O(\gd^3|\log\gd|^2)\, ,
\end{multline}
and this concludes the proof of the Theorem \ref{th:exp w} (the reversed inequality is evident, as already stated in the Introduction).

\end{subsection}

\medskip

\begin{subsection}{Proof of Corollary \lowercase{\ref{cor:w}}}\label{sec:proof cor}
We first justify the existence of an optimal path for each point $q_\gd(\varphi^\gd)$ with $\varphi^\gd\in \bbR/L_\gd\bbR$. We can not indeed apply directly Lemma \ref{lem:opt path}, since $q_\gd(\varphi^\gd)$ may not realize the minimum of the quasipotential on $\partial U^\gd$. But using similar compactness arguments as in the proof of Lemma \ref{lem:opt path}, it is sufficient to prove that 
\begin{equation}
\label{eq:cond existence opt path}
 W_\gd(A^\gd,q_\gd(\varphi^\gd))\, <\, \inf_{E\in \partial \tilde U^\gd}W_\gd(A^\gd,E)\, ,
\end{equation}
where $\tilde U^\gd$ is the whole tube introduced in Lemma \ref{lem:opt path}. Indeed this implies that a trajectory $Z$ linking $A^\gd$ to $q_\gd(\varphi^\gd)$ with $I_{\gd,\infty}(Z)$ sufficiently close to $W_\gd(A^\gd,q_\gd(\varphi^\gd))$ stays in $\tilde U^\gd$, and the Arzel\`a-Ascoli Theorem ensures the existence of an optimal path. For each $B^\gd\in \partial\tilde U^\gd$ satisfying
\begin{equation}
 W_\gd(A^\gd,B^\gd)\, =\,  \inf_{E\in \partial \tilde U^\gd}W_\gd(A^\gd,E)
\end{equation}
there exists an optimal path $Y^\gd$ staying in $\tilde U^\gd$ (see Remark \ref{rem:opt path tilde U}). Moreover this optimal path $Y^\gd$ satisfies the Lemmas \ref{lem:dot Y 1/2} to \ref{lem:dot h 2}. Indeed the proofs of these Lemmas do not rely on the position of the phase $\varphi_{B^\gd}$. So in particular $Y^\gd$ satisfies \eqref{ineq:I gd tau1 2}, and recalling \eqref{ineq:first bound Wred} and \eqref{eq:choice C0} this implies \eqref{eq:cond existence opt path}.
 
 Now an optimal path $Y^\gd$ associated to $q_\gd(\varphi^\gd)$ also satisfies the Lemmas \ref{lem:dot Y 1/2} to \ref{lem:dot h 2}, and this implies in particular that the associated 
 exit time $\tau^\gd_1$ is equal to $0$ (since $Y_0=q_\gd(\varphi^\gd)\in M^\gd$, and an optimal path that leaves the $\gd^2$-neighborhood of $M^\gd$ can not come back). To prove Corollary \ref{cor:w} it remains to do the same expansions as in the proof of Theorem \ref{th:exp w} (the expansion on the time interval $[\tau^\gd_1,0]$ is of course not needed here).
 
\end{subsection}

\end{section}

\section*{Acknowledgements}

This work is a part of my PhD Thesis. I would like to thank my advisor Giambattista Giacomin for having proposed this subject and for all his help and advices.
I would also like to thank Vincent Millot for very useful discussions.


\begin{thebibliography}{99}

\bibitem{cf:Berglund2010}
N. Berglund and B. Gentz, \emph{The Eyring-Kramers law for potentials with nonquadratic saddles}, Markov Processes Relat. Fields {\bf 16} (2010), 549-598.

\bibitem{cf:BEGK}
A. Bovier, M. Eckhoff, V. Gayrard, M. Klein, \emph{Metastability in reversible diffusion processes. I. Sharp asymptotics for capacities and exit times}, J. Eur. Math. Soc. (JEMS) {\bf 6} (2004), 399-424.


\bibitem{cf:Dayexpo}
M.V. Day, \emph{On the exponential exit law in the small parameter exit problem}, Stochastics {\bf 8} (1983), 297-323

\bibitem{cf:Daybound}
M. V. Day, \emph{Large deviations results for the exit problem with characteristic boundary}, J. Math. Anal. Appl. {\bf 147} (1990), 134-153. 

\bibitem{cf:DayDarden}
M.V. Day, T.A. Darden, \emph{Some regularity results on the Ventcel-Freidling quasi-potential function}, Appl Math Optim {\bf 13} (1985), 259-282.

\bibitem{cf:Eyring}
H. Eyring, \emph{The activated complex in chemical reactions}, Journal of Chemical Physics {\bf 3} (1935), 107-115.

\bibitem{cf:fen}
N. Fenichel, \emph{Persistence and smoothness of invariant manifolds for flows}, Indiana Univ. Math. J. {\bf 21} (1972), 193-226.

\bibitem{cf:FreidWentz}
M.I. Freidlin and A.D. Wentzell, \emph{Random Perturbations of dynamical systems}, Grundlehren der Mathematischen Wissenschaften Series, Springer Verlag, 1998.

\bibitem{cf:GOV}
A. Galves, E. Olivieri and M.E. Vares, \emph{Metastability for a class dynamical systems subject to small random perturbations}, The Annals of Probability {\bf 15} (1987), 1288-1305.

\bibitem{cf:GPPP}
G. Giacomin, K. Pakdaman, X. Pellegrin and C. Poquet, \emph{Transitions in  active rotator systems: invariant hyperbolic manifold  approach}, SIAM J. Math. Anal. {\bf 44}
(2012), 4165-4194.


\bibitem{cf:GTNE}
D.S. Goldobin, J. Teramae, H. Nakao, G.B. Ermentrout, \emph{Dynamics of limit-cycle oscillators subject to general noise}, Phys. Rev. Lett. {\bf 105} (2010), 154101.

\bibitem{cf:Wig}
G. Haller and I. Mezic and S. Wiggins,  \emph{Normally hyperbolic invariant manifolds in dynamical systems}, Applied Mathematical Sciences, Springer, 1994.

\bibitem{cf:HPS}
M.W. Hirsch and C.C. Pugh and M. Shub, \emph{Invariant manifolds}, Lecture Notes in Mathmatics {\bf 583}, Springer-Verlag, New York, 1977.

\bibitem{cf:Kur}
Y. Kuramoto, \emph{Chemical oscillations, waves, and turbulence}, Dover Books on Chemistry Series, Dover Publications, 2003.

\bibitem{cf:Kramer}
H.A. Kramer, \emph{Brownian motion in a field of force and the diffusion model of chemical reactions}, Physica {\bf 7} (1940), 284-304.

\bibitem{cf:MS}
R.S. Maier, D. L. Stein, \emph{A scaling theory of bifurcations in the symmetric weak-noise escape problem}, J. Statist. Phys. {\bf 83} (1996), 291-357.

\bibitem{cf:MOS}
F. Martinelli, E. Olivieri, E. Scoppola, \emph{Small random perturbations of finite and infinite dimensional dynamical systems: unpredictability of exit times}, J. Statist. Phys. {\bf 55} (1989), 477-504

\bibitem{cf:OV}
E. Olivieri and M.E. Vares, \emph{Large deviations and metastability}, Encyclopedia of Mathematics and its Applications, Cambridge University Press, 2005.

\bibitem{cf:SellYou}
G.R. Sell and Y. You, \emph{Dynamics of evolutionary equations}, Applied Mathematical Sciences {\bf 143}, Springer, 2002.

\bibitem{cf:TNE}
J. Teramae, H. Nakao, G.B. Ermentrout, \emph{Stochastic phase reduction for a general class of noisy limit cycle oscillators}, Phys. Rev. Lett. {\bf 102} (2009), 194102.

\bibitem{cf:YA}
K. Yushimura, K. Arai, \emph{Phase reduction of stochastic limit cycle oscillators}, Phys. Rev. Lett. {\bf 101} (2008), 154101.

\end{thebibliography}
\end{document}